\documentclass[a4paper, twoside, notitlepage, 12pt]{amsart} 
\usepackage{amsmath,amsfonts,amssymb,amsthm,mathrsfs} 

\usepackage[matrix,arrow]{xy}

\usepackage{verbatim}
\usepackage[latin1]{inputenc}   

\title[Recovering the good component of the Hilbert scheme]{Recovering the good component of the Hilbert scheme}
\author{Torsten Ekedahl}
\address{%
Department of Mathematics \\ 
Stockholm University \\ 
Stockholm \\ 
Sweden}
\email{teke@math.su.se}

\author{Roy Skjelnes}
\address{%
Department of Mathematics \\ 
KTH \\ 
Stockholm \\ 
Sweden}
\email{skjelnes@math.kth.se}


\DeclareMathOperator{\id}{id}

\DeclareMathOperator{\Spec}{Spec}
\DeclareMathOperator{\Proj}{Proj}

\newcommand\ili{\mathop{\hbox{\vtop{\setbox0=\hbox{lim}%
\dimen0\wd0\box0\vskip2pt
\nointerlineskip\hbox to \dimen0{\leftarrowfill}}}}\nolimits}
\newcommand{\ra}{\longrightarrow}

\newcommand{\eqbeg}{\begin{equation}}
\newcommand{\eqend}{\end{equation}}
\newcommand{\arbeg}[1]{\begin{array}{#1}}
\newcommand{\arend}{\end{array}}

\newcommand{\p}{\mathfrak{p}}

\newcommand{\Sgot}{\mathfrak{S}}
\newcommand{\A}{\mathfrak{a}}

\newcommand{\R}{\mathfrak{R}}

\newcommand{\fpr}{\ifmmode\mathrm{FPR}\else\textrm{FPR}\fi}

\newcommand{\CalB}{{\mathscr{B}}}
\newcommand{\CalC}{{\mathscr{C}}}
\newcommand{\CalI}{{\mathscr{I}}}

\newcommand{\CalO}{{\mathscr{O}}}
\newcommand{\CalE}{{\mathscr{E}}}

\newcommand{\CalH}{{\mathscr{H}}}

\newcommand{\T}{\operatorname{T}}
\newcommand{\TS}{\operatorname{TS}}

\newtheorem{thm}[subsection]{Theorem}
\newtheorem{lem}[subsection]{Lemma}
\newtheorem{cor}[subsection]{Corollary}
\newtheorem{prop}[subsection]{Proposition}
\newtheorem{deflem}[subsection]{Definition-Lemma}

\theoremstyle{definition}
\newtheorem{defn}[subsection]{Definition}

\theoremstyle{remark}
\newtheorem{rem}[subsection]{Remark}

\numberwithin{equation}{subsection}

\newcommand{\Hilb}{\mathrm{Hilb}}

\newcommand\good{\mathrm{G}^n_{X/S}}
\newcommand\n{\operatorname{n}}
\newcommand\bn{\operatorname{b}}

\newcommand\sHom{{\CalH\kern -0.25ex{ \mathit o\mathit m\/ }}}
\newcommand\Ext{{\mathrm{Ext}}}
\newcommand\sExt{{\CalE\kern -0.25ex{ \mathit x\mathit t\/ }}}

\begin{document}
\pagenumbering{arabic}

\begin{abstract} 

We give an explicit construction, for a flat map $X \to S$ of algebraic spaces,
of an ideal in the $n$'th symmetric product of $X$ over $S$. Blowing up this
ideal is then shown to be isomorphic to the schematic closure in the Hilbert
scheme of length $n$ subschemes of the locus of $n$ distinct points. This
generalises Haiman's corresponding result (\cite{haiman_catalan}) for the affine
complex plane. However, our construction of the ideal is very different from
that of Haiman, using the formalism of divided powers rather than representation
theory.

In the non-flat case we obtain a similar result by replacing the $n$'th
symmetric product by the $n$'th divided power product.
\end{abstract}

\maketitle

The Hilbert scheme, $\operatorname{Hilb}^n_{X/S}$, of length $n$ subschemes of a
scheme $X$ over some $S$ is in general not smooth even if $X \to S$ itself is
smooth. Even worse, it may not even be (relatively) irreducible. In the case of
the affine plane over the complex numbers (where the Hilbert scheme is smooth
and irreducible) Haiman (cf., \cite{haiman_catalan}) realised the Hilbert scheme
as the blow-up of a very specific ideal of the $n$'th symmetric product of the
affine plane.  It is the purpose of this article to generalise Haiman's
construction. As the Hilbert scheme in general is not irreducible while the
symmetric product is (for a smooth geometrically irreducible scheme over a field
say) it does not seem reasonable to hope to obtain a Haiman like description of
all of $\operatorname{Hilb}^n_{X/S}$ and indeed we will only get a description
of the schematic closure of the open subscheme of $n$ distinct points.  With
this modification we get a general result which seems very close to that of
Haiman. The main difference from the arguments of Haiman is that we need to
define the ideal that we want to blow up in a general situation and Haiman's
construction seems to be too closely tied to the $2$-dimensional affine space in
characteristic zero.

As a bonus we get that our constructions work very generally.  We have thus
tried to present our results in a generality that should cover reasonable
applications (encouragement from the referee has made us make it more general
than we did in a previous version of this article).

There are some rather immediate consequences of this generality. The first one
is that we have to work with algebraic spaces instead of schemes as otherwise
the Hilbert scheme (as well as the symmetric product) may not exist.  A second
consequence is that we find ourselves in a situation where existing references
do not ensure the existence of $\operatorname{Hilb}^n_{X/S}$ and we give an
existence proof in the generality required by us (which is a rather easy
patching argument to reduce it to known cases).

It turns out that the key to constructing the ideal to blow up is to use the
formalism of divided powers. Recall that if $A$ is a commutative ring and $F$ a
flat $A$-algebra, then the subring of $\Sgot_n$-invariants of $F^{\otimes_{A}n}$
is isomorphic to the $n$'th divided power algebra $\Gamma^n_{A}(F)$ (through the
map that takes $\gamma^n(r)$ to $r^{\otimes n}$). 

Using the fact that $\Gamma^n(F)$ is the degree $n$ component of the divided
power algebra $\Gamma^*(F)$ we can define an ideal in $\Gamma^n(F)$ (this graded
component of the divided power algebra becomes an algebra using the
multiplication of $F$) which is our candidate to be blown up. Note that in the
definition of this ideal we are using in an essential way the multiplication in
the divided power algebra $\Gamma^*(F)$ forcing us to carefully distinguish
between the multiplication in this graded algebra and the multiplication of its
graded component $\Gamma^n(F)$ induced by the multiplication on $F$. On the
upside it is exactly this interplay that allows us to define, in a generality
outside of Haiman's case, the ideal. Furthermore, the excellent formal
properties of $\Gamma^n(F)$ allows us to \emph{define} an analogue of the
symmetric product of $\Spec(F) \to\Spec(A)$ as $\Spec(\Gamma^n(F)) \to \Spec(A)$ in
the case when $A \to F$ is not flat. This makes our arguments go through without
problems in the case when $\Spec(F) \to \Spec(A)$ is not necessarily flat. (We
also need to extend the construction of $\Spec(\Gamma^n(F))$ to the non-affine
case; the gluing argument needed to make this extension uses results of David
Rydh, \cite{rydh_fams}.)

In more detail this paper has the following structure: We start with some
preliminaries on divided powers and recall of the Grothendieck-Deligne norm
map. The main technical result is to be found in Sections 5 and 6. There we
first find a (local) formula for the multiplication of the tautological rank
$n$-algebra over the configuration space of $n$ distinct points of $X$. We then
note that this formula makes sense over the blow-up of a certain ideal in the
full symmetric product. This gives us a family of length $n$ subschemes of $X$
over this blow-up and hence a map of it to the Hilbert scheme. Once having
constructed it, it is quite easy to show that it gives an isomorphism of the
blow-up to the schematic closure of the subspace of $n$ distinct points of the
Hilbert scheme. The proof first does this in the case $X \to S$ is affine and
then discusses the patching (and limit arguments) needed to extend it to the
more general case.

We finish by tying some loose ends. First we generalise the result of Fogarty on
the smoothness of $\operatorname{Hilb}^n_{X/S}$ for $X \to S$ smooth of relative
dimension $2$ removing the conditions on the base $S$ needed by
Fogarty. Finally, we discuss how one can, under suitable conditions, embed the
blow-up in a Grassmannian as Haiman does.

\subsection{Acknowledgments} 
First we would like to give our thanks to the referee of an earlier version of
this article who insisted on putting our results in the most general context and
who also sketched how to make such an extension as well as pointing out
inaccuracies in the article.

We have received a lot of inspiration with conversations with our colleagues.
In particular the second author would like to acknowledge input from Steve
Kleiman, Ezra Miller and Dan Laksov.

\section{Divided powers and norm}

In this section we first recall some  properties for the ring of
divided powers. The standard reference is Roby \cite{roby_lois_pol_mod}
and
\cite{roby_lois_pol_mult}, but see also \cite{bourbaki_polyetfractions} and
\cite{ferrand_norme}. Algebras in this note are commutative.

\subsection{The ring of divided powers}\label{subsec1.1} Let $A$ be a commutative ring
and $M$ an $A$-module. The ring of divided powers $\Gamma_AM$ is
constructed as follows. We consider the polynomial ring over
$A[\gamma^{n}(x)]_{(n,x)\in \bold{N}\times M}$, where the variables
$\gamma^n(x)$ are indexed by the set $\bold{N}\times M$, where
$\bold{N}$ is the set of non-negative integers. Then the ring
$\Gamma_AM$ is obtained by dividing out the polynomial ring by the
following relations 
\begin{align}
 & \gamma^0(x)-1 \label{111}\\
 & \gamma^n(\lambda x)-\lambda^n \gamma^n(x) \label{112}\\
 & \gamma^n(x+y)-\sum_{j=0}^n\gamma^j(x)\gamma^{n-j}(y) \label{113}\\
 & \gamma^n(x)\gamma^m(x)-\binom{n+m}{n}\gamma^{n+m}(x) \label{114}
\end{align}
for all integers $m,n\in \bold{N}$, all $x,y\in M$, and all $\lambda
\in A$.  The residue class of the variable $\gamma^n(x)$ in
$\Gamma_AM$ we denote by $\gamma^n_M(x)$, or simply $\gamma^n(x)$ if
no confusion is likely to occur. The ring $\Gamma_AM $ is graded
where $\gamma^n(x)$ has degree $n$, and with respect to this grading
we write $\Gamma_AM=\oplus_{n\geq 0}\Gamma^n_AM$.

\subsection{Polynomial laws } Let $A$ be a ring, and let $M$ and $N$ be
two fixed $A$-modules. Assume that we for each $A$-algebra $B$ have a
map {\it of sets} $g_{B}
\colon M\otimes_AB \ra N\otimes_AB$ such that for any $A$-algebra homomorphism $u \colon B\ra
B'$ the following diagram is commutative
$$ \xymatrix{ M\otimes_AB \ar[d] \ar[r]^{g_B} & N\otimes_AB\ar[d] \\
 M\otimes_AB' \ar[r]^{g_{B'}} & N\otimes_AB', }
$$
where the vertical maps are the canonical homomorphisms. Such a collection of maps is called a polynomial law from $M$ to $N$,
and we denote the polynomial law with $\{g\} \colon M\ra N$.

\begin{defn}[Norms] Let $A$ be a ring, $M$ and $N$ two $A$-modules.
\renewcommand{\labelenumi}{(\theenumi)}
\begin{enumerate}
\item A polynomial law $\{g\}\colon M \ra N$ is {\it homogeneous of degree} $n$ if
  for any $A$-algebra $B$ we have that $g_{B}(bx)=b^ng_B(x)$, for any
  $x\in M\otimes_AB$ and any $b\in B$.
\item A  polynomial law $\{g\}\colon F\ra E$ between two $A$-algebras $F$
  and $E$, is {\it multiplicative} if
  $g_B(xy)=g_B(x)g_B(y)$ for any $x$ and $y$ in $F\otimes_AB$, for any $A$-algebra $B$. Furthermore, we require that $g_B(1)=1$.
\end{enumerate}
A norm (of degree $n$) from an $A$-algebra $F$ to an $A$-algebra $E$ is a homogeneous multiplicative polynomial law of degree $n$.
\end{defn}

\subsection{Universal norms}\label{sec1.4} Let $n$ be
a non-negative integer. For any $A$-algebra $B$ we have that
$\Gamma^n_A(M) \otimes_AB$ is canonically identified with
$\Gamma^n_B(M\otimes_AB)$. It follows that we have a polynomial law
$\{\gamma^n\} \colon M \ra \Gamma^n_AM$ and by (\ref{112}) the law is
homogeneous of degree $n$. The polynomial law $\{\gamma^n\} \colon M\ra
\Gamma^n_AM$ is {\it universal} in the sense that the assignment
$u\mapsto \{u\circ \gamma^n\}$ gives a bijection between the
$A$-module homomorphisms $u\colon \Gamma^n_AM \ra N$ and the set of
polynomial laws of degree $n$ from $M$ to $N$.

Furthermore, if $F$ is an $A$-algebra then $\Gamma^n_AF$ is an
$A$-algebra and then the polynomial law $\{\gamma^n \} \colon F\ra
\Gamma^n_AF$ is the universal {\it norm} of degree $n$ (\cite[Thm.\ p.\ 871]{roby_lois_pol_mult}, \cite[2.4.2, p.\ 11]{ferrand_norme} ). ``Universal''
here means in the sense as described above, but for $A$-algebra
homomorphisms from $\Gamma^n_AF$.

\subsection{The different products} The product structure on
$\Gamma_AF$ we refer to as the external structure. We will denote
the external product with $*$ in order to distinguish the external
product from the product structure on each graded component
$\Gamma^n_A F$ defined in the previous section. (Note that our
convention is the reverse of the one used in \cite{ferrand_norme}.)
 
\subsection{} Let $p$ a fixed positive integer, and for
each $i=1, \ldots, p$ we assume that we have a sequence of
non-negative integers $\{a_{i,j}\}$ $(1\leq j\leq q_i)$ such that
$\sum_{j=1}^{q_i}a_{i,j}=n$. Let $I=[1,
\ldots, q_1]\times \cdots \times [1,\ldots , q_p]$, and let $\CalB\{a_{i,j}\}\subset \times_{I}\bold{N}$ be the subset of integers $b=\{b_{i_1, \ldots,
  i_p}\}_{i_1,i_2, \ldots, i_p \in I}$ such that their sum
  $|b|:=\sum_{i_1,\ldots, i_p \in I}b_{i_1,\ldots, i_p}=n$ and such
  that
\eqbeg \sum b_{i_1, \ldots, i_{r-1},s,i_{r+1}, \ldots
,i_p}=a_{r,s},\label{171}
\eqend for all $r,s$.

\begin{prop}\label{prop1.7} Let $q_1, \ldots, q_p$ be $p$ positive integers, and let $\{a_{i,j}\}$ be any sequence of non-negative integers such that
$\sum_{j=1}^{q_i}a_{i,j}=n$ for each $i=1, \ldots, p$. Then we have for any
elements $\{x_{i,j}\}$  in the
$A$-algebra $F$ that in $\Gamma^n_AF$ the following identity holds
$$ \big( \prod_{j=1}^{q_1}*\gamma^{a_{1,j}}(x_{1,j})\big)\cdots \big(
\prod_{j=1}^{q_p}* \gamma^{a_{p,j}}(x_{p,j})\big)$$
$$=\sum_{\CalB\{a_{i,j}\}} \prod_{i_1, \ldots, i_p \in I}
    *\gamma^{b_{i_1, \ldots, i_p}}(x_{1,i_1}x_{2,i_2}\cdots
    x_{p,i_p}).$$
\end{prop}   

\begin{proof} For each $i=1, \ldots ,p$ we consider the linear polynomials 
$$ L_i=x_{i,1}T_{i,1}+x_{i,2}T_{i,2} +\dots +x_{i,q_i}T_{i,q_i}$$ in
the variables $\{T_{i,j}\} (1\leq j\leq q_i)$. Let us denote by $A[T]$
the polynomial ring in the variables $\{T_{i,j}\}$ over $A$, and let
$F[T]=A[T]\otimes_AF$. By iterating the formula (\ref{113}) and using
the fact that the variables are scalars over $A[T]$ combined with
(\ref{112}) we achieve the following expression 
$$ \gamma^n(L_i)=\sum_{|\alpha_{i}|=n}\big(\prod_{j=1}^{q_i}*\gamma^{\alpha_{i,j}}(x_{i,j})\big)T_{i,1}^{\alpha_{i,1}}\cdots T_{i,q_i}^{\alpha_{i,q_i}} \quad \in \Gamma^n_{A[T]}F[T],$$
where we have abbreviated $\alpha_{i,1}+\cdots +
\alpha_{i,q_{i}}=|\alpha_i|$. We therefore get 
\eqbeg
\gamma^n(L_1)\cdots
\gamma^n(L_p)=\prod_{i=1}^p\big(\sum_{|\alpha_{i}|=n}\big(\prod_{j=1}^{q_i}*\gamma^{\alpha_{i,j}}(x_{i,j})\big)T_{i,1}^{\alpha_{i,1}}\cdots
T_{i,p_i}^{\alpha_{i,q_i}}\big). \label{172}
\eqend

On the other hand we have   $\gamma^n(L_1)\cdots
\gamma^n(L_p)=\gamma^n(L_1\cdots L_p)$, thus the identity (\ref{172}) above also equals
\eqbeg \gamma^n(L_1\cdots L_p) =\gamma^n(\sum_{i_1, \ldots, i_p \in
I}x_{1,i_1}x_{2,i_2}\cdots x_{p,i_p}T_{1,i_1}T_{2,i_2}\cdots
T_{p,i_p}), \label{173}
\eqend
with $I=[1,\ldots, q_1]\times \cdots \times [1, \ldots ,q_p]$. We then
finally iterate the right side of the expression (\ref{173}) by the
formulas (\ref{112}) and (\ref{113}), and obtain that (\ref{173}) can be written as
\eqbeg \sum_{|b|=n} \big(\prod_{i_1, \ldots ,i_p \in I}*\gamma^{b_{i_1,
    \ldots ,i_p}}(x_{1,i_1} x_{2,i_2}\cdots
    x_{p,i_p})T_{1,i_1}^{b_{1,i_1, \ldots ,p,i_p}}\cdots
    T_{p,i_p}^{b_{i_1, \ldots ,i_p}}\big).\label{174}
\eqend
The sum is to be taken over all integers $b=\{b_{i_1, \ldots ,i_p}\}$
in $\times_I \bold{N}$ such that $|b|=n$. We now compare the
coefficients of the polynomial (\ref{172}) with the polynomial in
(\ref{174}). The coefficient of
$$ T_{1,1}^{a_{1,1}}\cdots T_{1,q_1}^{a_{1,q_1}}\cdots T^{a_p,1}_{p,1}
\cdots T_{p,q_p}^{a_{p,q_p}}$$ 
yields the result.  
\end{proof}

\begin{rem} The proposition is a generalization of the case when $p=2$ that can be found in (\cite[Formula 2.4.2]{ferrand_norme}). The trick we used in the proof is to add variables and then iterate the defining equations for the ring of divided powers. The technique was communicated to us by  Dan Laksov. 
\end{rem}

As an example of statement in the proposition we give here an identity that we will use later. 
\begin{lem}\label{gammaf^n} Let $x_1, \ldots ,x_n$ and $f$ be elements in an $A$-algebra $F$. Then we have that $\gamma^1(x_1f^n)*\gamma^1(x_2) * \cdots * \gamma^1(x_n)$ equals
$$ \sum _{c=1}^n(-1)^{c+1} (\gamma^{c}(f) * \gamma^{n-c}(1)) \cdot (\gamma^1(x_1f^{n-c})* \gamma^1(x_2) * \cdots * \gamma^1(x_n)).
$$
\end{lem}
\begin{proof} For each $c=1, \ldots, n$ the above proposition identifies the product $\gamma^c(z_{1,1})*\gamma^{n-c}(z_{1,2})\cdot \gamma^1(z_{2,1}) *\cdots * \gamma^1(z_{2,n})$ with 
\eqbeg \label{expansion} \sum_{\CalB \{a_{i,j}^c\}} \prod * \gamma^{b_{i_1,i_2}}(z_{1,i_1}z_{2,i_2}),
\eqend
where the sequence $\{a_{i,j}^c\}$ is $a_{1,1}^c=c, a_{2,1}^c=n-c$, and $a_{2,j}^c=1$ for $j=1, \ldots ,n$. Let $J=[1,2]\times [2, \ldots, n]$, and define $\CalC^c \subset \times_J \{0,1\}$ as
$$ \CalC^c =\{ \{b_{j_1,j_2} \} \mid \sum_{j=2}^nb_{1,j}=c, \sum_{j=2}^nb_{2,j}=n-c-1, b_{1,j}+b_{2,j}=1\}.$$
It is then clear  that the indexing set $\CalB \{a_{i,j}^c\}$ is the disjoint  union
$$ \CalB \{a_{i,j}^c\} =\{b_{1,1}=1, b_{2,1}=0\}\times \CalC^{c-1} \sqcup \{ b_{1,1}=0, b_{2,1}=1\}\times \CalC^c.$$ 
We then have that the expression (\ref{expansion}) can be written as
$$ \sum_{\CalC^{c-1}}\gamma^1(z_{1,1} z_{2,1})*\prod \gamma^{b_{j_1,j_2}}(z_{1,j_1}z_{2,j_2}) + \sum_{\CalC^{c}}\gamma^1(z_{1,2}z_{2,1})*\prod \gamma^{b_{j_1,j_2}}(z_{1,j_1}z_{2,j_2}).$$
Assume now that we have, for each $c=1, \ldots, n$, a collection of elements $\{z_{i,j}^c\}$ with $(i,j)\in [1,2]\times [1, \ldots , n]$ such that $z_{1,2}^cz_{2,1}^c=z_{1,1}^{c+1}z_{2,1}^{c+1}$ for $c=1, \ldots, n-1$,  and that $z_{1,j_1}^{c}z_{2,j_2}^c=z_{1,j_1}^{c+1}z^{c+1}_{2,j_2}$ for $(j_1,j_2) \in [1,2]\times [2,\ldots ,n]$. Using the splitting of $\CalB \{a_{i,j}^c\}$ described above we have that the alternating sum
$$\sum_{c=1}^n (-1)^{c+1}\sum_{\CalB \{a_{i,j}^c\}} \prod * \gamma^{b_{i_1,i_2}}(z^c_{1,i_1}z^c_{2,i_2})$$ 
is a telescoping sum. We note that $\CalC^n $ is the empty set, and that $\CalC^0$ is the singleton set $\CalC^0 =\{ \{b_{1,j}=0, b_{2,j}=1\} \mid j=2, \ldots, n\}$. Consequently the telescoping sum collapses to
$$\gamma^{1}(z^1_{1,1}z^1_{2,1})*\gamma^{1}(z^1_{1,2}z^1_{2,2})* \cdots * \gamma^1(z^1_{1,2}z^1_{2,n}).$$
The lemma is the special situation with $z_{1,1}^c=f, z_{1,2}^c=1, z_{2,1}^c=x_1f^{n-c}$ and $z^c_{2,j}=x_j$ for $j=2, \ldots , n$, $c=1, \ldots, n$.
\end{proof}

\subsection{The canonical homomorphism} An important norm is the following. Let $E$ be an $A$-algebra that is free of finite rank $n$ as an
$A$-module. For any $A$-algebra $B$ we have the determinant map $d_B \colon
E\otimes_AB \ra B$ sending $x\in E\otimes_AB$ to the determinant of
the $B$-linear endomorphism $e\mapsto ex$ on $E\otimes_AB$. It is
clear that the determinant maps give a multiplicative polynomial law $\{d\} \colon E \ra
A$, homogeneous of degree $n=\text{rank}_AE$. By
the universal properties (\ref{sec1.4}) of $\Gamma^n_AE$  we then have an
$A$-algebra homomorphism 
\eqbeg \sigma_E \colon \Gamma^n_AE \ra A, \eqend
such that $\sigma_E(\gamma^n(x))=\det (e\mapsto ex)$ for all $x\in
E$. In fact we get a homomorphism $\sigma_E \colon \Gamma^n_AE
\ra A$ even if $E$ is only locally free of finite rank $n$. We call
$\sigma_E$ the canonical homomorphism (\cite[Section 6.3, p.180]{sga4_deligne_coh_supp_prop}, \cite[Section 1.4, p.13]{iversen_lineardeterminants}).

\begin{prop}\label{prop1.11} Let $E$ be an $A$-algebra such that $E$ is
free of finite rank $n$ as an $A$-module. For any element $x\in E$ the
characteristic polynomial $\det(\Lambda -x) \in A[\Lambda]$ of the endomorphism $e\mapsto ex$ on $E$ is
$ \Lambda^n +\sum_{j=1}^n
(-1)^j\Lambda^{n-j}\sigma_E(\gamma^{j}(x)*\gamma^{n-j}(1))$.
In particular we have 
$$\operatorname{Trace}(e\mapsto ex) = \sigma_E
(\gamma^{1}(x)*\gamma^{n-1}(1)).$$
\end{prop}

\begin{proof} Let $\Lambda $ be an independent variable over $A$, and
write $E[\Lambda]=E\otimes_A A[\Lambda ]$. By the defining property of
the canonical homomorphism $\sigma_{E[\Lambda]}$ we have that the
characteristic polynomial $\det (\Lambda
-x)=\sigma_{E[\Lambda]}(\gamma^n(\Lambda-x))$.  We now use the
defining relations (\ref{112}) and (\ref{113}) in the
$A[\Lambda]$-algebra $\Gamma^n_{A[\Lambda]}E[\Lambda]$ and obtain

\begin{align}
\gamma^n(\Lambda-x)&=\sum_{j=0}^n(-1)^j\gamma^j(x)*\gamma^{n-j}(\Lambda)\notag \\
&=\sum_{j=0}^n(-1)^j\Lambda^{n-j}\gamma^j(x)*\gamma^{n-j}(1).\notag
\end{align} 

We have that $\Gamma^n_A(R)\otimes_A B =\Gamma^n_B(R\otimes_AB)$ and that
$\sigma_{E[\Lambda]}=\sigma_E\otimes\id_{A[\lambda]}$. Consequently
$\sigma_{E[\Lambda]}$ acts trivially on the variable $\Lambda$ and that the
action otherwise is as $\sigma_E$. Thus we obtain that
$\sigma_E(\gamma^j(x)*\gamma^{n-j}(1))$ in $A$ is the $j$'th coefficient of the
characteristic polynomial of $e\mapsto ex$ which proves the claim.
\end{proof}

\section{Discriminant and ideal of norms}

In this section we define the important ideal of norms and show their
connection with discriminants.

\begin{defn} Let $F$ be an
$A$-algebra. For any $2n$-elements $x=x_1, \ldots ,x_n$ and $y=y_1, \ldots ,y_n$ in $F$ we define $\delta (x,y)\in \Gamma^n_AF$ as
the element
$$ \delta (x,y):=\det {*} [\gamma^1 (x_iy_j)]_{1\leq i,j\leq n}=\sum_{\sigma \in \Sgot_n}(-1)^{|\sigma|}\gamma^1(x_1y_{\sigma(1)})* \cdots * \gamma^1(x_ny_{\sigma(n)}).$$
\end{defn}
 
\begin{rem}Note that for each element $z\in F$ the element $\gamma^1(z)$ is in
$\Gamma^1_AF=F$, but the product $\gamma^1(z_1)*\dots *\gamma^1(z_n)$
is in $\Gamma^n_AF$. 
\end{rem}

\begin{lem}\label{lemma2.2} Let $x_1, \ldots ,x_n$  and $y_1, \ldots, y_n$ be $2n$-elements in $F$. Then we have 
$$ \det *[\gamma^1(x_iy_j)]_{1\leq i,j \leq n}=\det[\gamma^{1}(x_iy_j)*\gamma^{n-1}(1)]_{1\leq i,j \leq n}.$$
\end{lem}

\begin{proof} We have that the determinant of $[\gamma^1(x_iy_j)*\gamma^{n-1}(1)]$ is 
\eqbeg\sum_{\sigma \in \Sgot_n}(-1)^{|\sigma|}\big(\gamma^{1}(x_1y_{\sigma(1)})*\gamma^{n-1}(1)\big) \cdots \big(\gamma^1(x_ny_{\sigma(n)})*\gamma^{n-1}(1)\big).\label{221}
\eqend
We will now use Proposition (\ref{prop1.7}) to expand the product in (\ref{221}). For that purpose we denote $I=[1,2]^n$, $a_{i,1}=1$ and $a_{i,2}=n-1$, and  we write $X^{\sigma}_{i,1}=x_iy_{\sigma (i)}$ and $X^{\sigma }_{i,2}=1$ for all $i=1, \ldots, n$. Using Proposition (\ref{prop1.7}), we then rewrite  the expression (\ref{221}) as
\eqbeg
\sum_{\sigma \in \Sgot_n} (-1)^{|\sigma|}\sum_{b\in
\CalB\{a_{i,j}\}}\prod_{i_1, \ldots ,i_n \in I}*
\gamma^{b_{i_1,\ldots ,i_n}}(X^{\sigma}_{1,i_1}\dots
X^{\sigma}_{n,i_n}).\label{222}
\eqend

 Consider now the element $\bold{b}
\in \times_I{\bold N}$ that is zero on all components except the
components indexed by $i_1, \ldots ,i_n$ where all indices but one
$i_r=1$ (and thus the other indices are $=2$), and let the values of
$\bold{b}$ at those components all equal 1. That is
$$ b_{1,2\ldots ,2}=b_{2,1,2, \ldots, 2}= \dots =b_{2,\ldots
,2,1}=1.$$ We then have that $\bold{b}$ satisfies the equations
(\ref{171}), hence $\bold{b} \in \CalB$. We see that $\delta(x,y)$
occurs as a summand in (\ref{222}) corresponding to having $\bold{b}$
fixed. We can therefore write (\ref{222}) as
\eqbeg
\delta(x,y)+\sum_{\sigma \in \Sgot_n} (-1)^{|\sigma|}\sum_{b\in
\CalB\{a_{i,j}\}\setminus \bold{b}}\prod_{i_1, \ldots ,i_n \in I}*
\gamma^{b_{i_1,\ldots ,i_n}}(X^{\sigma}_{1,i_1}\dots
X^{\sigma}_{n,i_n}) \label{223}
\eqend

We need to show that the right hand
side of (\ref{223}) is zero. By using the equations described by
(\ref{171}) one gets that if $b_{i_1, \ldots , i_n}$ is a non-zero
integer with at least one index $i_r=1$, then because $a_{r,1}=1$ we
must have the value $b_{i_1, \ldots , i_r=1, \ldots ,
i_n}=1$. Furthermore it follows easily that if $b\in
\CalB\{a_{i,j}\}$ where one component $b_{i_1, \ldots ,i_n}=1$ with
only one index $i_r=1$, then $b=\bold{b}$ defined above. Thus a
non-zero component of an element $b\in \CalB\{a_{i,j}\}\setminus
\bold{b}$ has either all indices $i_1=\dots =i_n=2$, or at least two
indices $i_r=i_{r'}=1$.

The first alternative with $i_1=\dots =i_n=2$ gives that the only
non-zero component of $b$ is $b_{2,\ldots ,2}$, and that alternative
is impossible as it should equal $n-1$ by (\ref{171}), but also it
should equal $|b|=n$, being an element of $\CalB\{a_{i,j}\}$.

Let us now rule out the other alternative, with at least two indices
$i_r=i_{r'}=1$ of an element $b\in \CalB\{a_{i,j}\}$.  The equations (\ref{171})
then imply that all the other non-zero components must have the
indices $i_r=i_{r'}=2$. From that we deduce that the product
$$  \prod_{i_1, \ldots, i_n \in I}*\gamma^{b_{i_1, \ldots ,i_n}}(X^{\sigma}_{1,i_1}\dots X^{\sigma}_{n,i_n})$$
is invariant under the action of $\epsilon$, the operation that permutes the two factors $r$ and $r'$. Because in one component of the product (corresponding to $i_r=i_{r'}=1$) we have
$$X^{\sigma}_{1,i_r}X^{\sigma}_{1,i_{r'}}=x_{i,r}y_{\sigma(i_r)}x_{i,r'}y_{\sigma(i_r')}=X^{\epsilon \sigma}_{1, i_r}X^{\epsilon \sigma}_{1, i_r'}.$$
In the other components of the product we have  $X^{\sigma}_{2,i_r}=X^{\sigma}_{2,i_r'}=1$, clearly invariant under permutation. As $\epsilon $ has sign -1 it is clear that the above sum is annihilated in the determinant expression. As the above sum was arbitrary we have that the right side of (\ref{223}) is zero. 
\end{proof}

\begin{lem}\label{product} Let $x=x_1, \ldots ,x_n$ and $y=y_1, \ldots, y_n$ be $2n$ elements in an $A$-algebra $F$. We have $\delta(x,y)^2=\delta(x,x)\delta(y,y)$.
\end{lem}

\begin{proof} By definition $\delta(x,y)^2$ is the sum
$$ \sum _{\sigma, \tau \in \Sgot_n}(-1)^{|\sigma \tau|} \big( \prod_{i=1}^n * \gamma^1(x_iy_{\sigma (i)})\big) \cdot  \big( \prod_{j=1}^n *\gamma^1(x_jy_{\tau (j)})\big).$$
We expand each summand  using Proposition (\ref{prop1.7}) which gives that $\delta (x,y)^2$ can be written as
$$ \sum_{\sigma, \tau \in \Sgot_n}(-1)^{|\sigma \tau|} \sum_{\rho \in \Sgot_n}\prod _{i=1}^n * \gamma^1(x_iy_{\sigma(i)}x_{\rho(i)}y_{\tau(\rho (i))}).$$
On the other hand we have that $\delta(x,x)\delta(y,y)$ is
$$ \sum_{\rho, \beta \in \Sgot_n}(-1)^{|\rho \beta|} \big( \prod_{i=1}^n *\gamma^1(x_ix_{\rho(i)}) \big) \cdot \big( \prod_{j=1}^n* \gamma^1(y_jy_{\beta (j)})\big).$$
Expanding the products as in Proposition (\ref{prop1.7}), and then substituting $\beta =\tau \rho \sigma^{-1}$, proves the lemma.
\end{proof}

\begin{lem}\label{lemma2.3} Let $x=x_1, \ldots ,x_n$ and $y=y_1, \ldots, y_n$ be $2n$ elements in an $A$-algebra $F$. Let $z_j=\sum_{i=1}^na_{i,j}x_i$ and $w_j=\sum_{i=1}^nb_{i,j}y_i$ be $A$-linear combinations of $x$ and $y$ for $j=1, \ldots, n$. Then we have
$$\delta (z, w)=\det (a_{i,j})\det(b_{i,j})\delta(x,y).$$
\end{lem}

\begin{proof} We clearly have the two matrix equations
\eqbeg
 Z=(a_{i,j})X \quad \text{ and }\quad W=(b_{i,j})Y,\label{231}
\eqend

where $X,Y,Z$ and $W$ are the $(n\times 1)$-matrices with entries $x_i,y_i,z_i$ and $w_i$ ($1\leq i\leq n$), respectively.

We have $ax=a\gamma^1(x)$ in $\Gamma^1_AF=F$ for any $a\in A$ and any $x\in F$. We therefore view (\ref{231}) as a matrix equation over $\Gamma^1_AF$. We then multiply the matrix $Z$ with the transpose of $W$ and obtain an $(n\times n)$-matrix
\eqbeg
 Z \cdot W^{tr}=(a_{i,j}) X\cdot Y^{tr} (b_{i,j}). \label{232}
\eqend
As we have $\gamma^1(z)\cdot \gamma^1(w)=\gamma^1(zw)$ in $\Gamma^1_AF$ we read (\ref{232}) as
\eqbeg
 (\gamma^{1}(z_iw_j)) =(a_{i,j}) (\gamma^1(x_iy_j)) (b_{i,j})^{tr}.\label{233}
\eqend
The $A$-algebra structure on the ring $\Gamma_AF$ is compatible with the $A$-algebra structure on each of its graded components $\Gamma^m_AF$, and we will therefore view (\ref{232}) as an equation of four matrices $Z\cdot W^{tr}, (a_{i,j}), X\cdot Y^{tr}$ and $(b_{i,j})$ over the $A$-algebra $\Gamma_AF$. Now, using the usual properties of the determinant we obtain the following identity
$$\det *(\gamma^1(z_iw_j))=(\det *(a_{i,j}))*(\det *(\gamma^1(x_iy_j)))*(\det *(b_{i,j})).$$
As the entries of $(a_{i,j})$ and the entries of $(b_{i,j})$ are in the ring $A=\Gamma^0_AF$ we get that $\det *(a_{i,j})=\det (a_{i,j})$, and similarly that $\det *(b_{i,j})=\det (b_{i,j})$. It then follows that the last displayed equation is what we wanted to prove. 
\end{proof}

\begin{defn}[The ideal of norms]\label{defn2.4} Let $n$ be a fixed integer, and let $V\subseteq F$ be an $A$-submodule of an $A$-algebra $F$. We define $I_V\subseteq \Gamma^n_AF$, the ideal of norms associated to $V$, as the ideal generated by
$$ \delta(x,y)=\det * [\gamma^{1}(x_iy_j)]_{1\leq i,j \leq n} \in \Gamma^n_AF$$
for any $2n$-elements $x=x_1, \ldots ,x_n$ and $y=y_1, \ldots ,y_n$ in $V\subseteq F$. 
\end{defn}

\begin{lem}\label{lemma2.5} Let $A \ra B$ be a homomorphism of rings, and let $V\subseteq F$ be an $A$-submodule of an $A$-algebra $F$. The extension of the ideal $I_V$ by the $A$-algebra homomorphism $ \Gamma^n_AF \ra \Gamma^n_A(F)\otimes_AB$ equals the ideal $I_{V_B}$; the ideal of norms associated to the $B$-submodule $\text{Im}(V\otimes_AB \ra F\otimes_AB)$.
\end{lem}

\begin{proof} Let $x=x_1, \ldots ,x_n$ and $y=y_1, \ldots, y_n$ be $2n$-elements in $V\subseteq F$. The element $\delta (x,y)$ in $\Gamma^n_AF$ is then mapped to $\delta (x,y)\otimes 1_B =\det *[\gamma^1_F(x_iy_j)\otimes 1_B]$
in $\Gamma^n_A(F)\otimes _AB$. We have that $\Gamma^n_A(F)\otimes _AB$ is canonically identified with $\Gamma^n_B(F\otimes_AB)$. 
Hence
$$\det *[\gamma^1_F(x_iy_j)\otimes 1_B]_{1\leq i,j \leq n}=\det *[\gamma^1_{F\otimes_AB}(x_iy_j \otimes 1_B)]_{1\leq i,j \leq n}.$$
We then have that the extension $I_V\Gamma^n_B(F\otimes_AB)$ is included in $I_{V_B}$. As the generators of the $B$-module $V_B$ are the images of generators of the $A$-module $V$, the inclusion $I_{V_B}\subseteq I_V\Gamma^n_B(F\otimes_AB)$ follows from Lemma (\ref{lemma2.3}). 
\end{proof}

\begin{lem}\label{modpolyring} Let $F=A[T_1, \ldots ,T_r]$ be the polynomial
ring in a finite set of variables, and let $V\subset F$ be the $A$-module
spanned by those monomials whose degree in each of the variables is less than
$n$. Then the ideals of norms associated to $V$ and $F$ are equal; that is
$I_V=I_F$. Furthermore, if $n!$ is invertible in $A$ then $I_W=I_F$, where
$W\subset F$ is the $A$-module spanned by monomials of degree less than $n$.
\end{lem}

\begin{proof} Given $x_1, \ldots ,x_n$ and $f$ in  $F$ we write $x(c)=x_1f^c, x_2, \ldots, x_n$. For any $y_1, \ldots, y_n$ we then obtain from the equality given in Lemma (\ref{gammaf^n}) that
$$ \delta (x(n),y)= \sum_{c=1}^n(-1)^{c+1}(\gamma^c(f) * \gamma^{n-c}(1)) \cdot \delta (x(n-c),y).$$
The first assertion of the lemma follows from the above equality. When $n!$ is invertible, the $n$'th powers of linear forms span the module of degree $n$ monomials, and the above equality then also yields the second assertion. 
\end{proof}

\subsection{Discriminant}  Let $E$ be an $A$-algebra
that is free of finite rank $n$ as an $A$-module. The trace of an $A$-linear
endomorphism of $E$ is an $A$-linear map $\operatorname{End}_A(E)\ra A$, that
composed with the natural map $E \ra \operatorname{End}_A(E)$ sends an element
$x\in E$ to $\operatorname{tr}(x)$; the trace of the endomorphism $e\mapsto
ex$ on $E$. There is an associated map $E \ra \operatorname{Hom}_A(E,A)$
taking $y\in E$ to the trace $\operatorname{tr}(xy)$, for any $x\in E$. 

The discriminant ideal
$D_{E/A}\subseteq A$ is defined (see e.g. \cite[p. 124]{altklei_gr_duality}) as the
ideal generated by the determinant of the associated map $E \ra
\operatorname{Hom}(E,A)$. 

\begin{prop}\label{prop2.7} Let $E$ be an $A$-algebra that is free of finite rank $n$ as an $A$-module. Then the extension of  $I_V$, the ideal of norms associated to $V=E$, by the canonical homomorphism $\sigma_E \colon \Gamma^n_AE \ra A$ is the discriminant ideal. In particular we have that the extension $\sigma_E(I_V)A=A$ if and only if $\Spec (E) \ra \Spec (A)$ is \'etale.
\end{prop}

\begin{proof} By Lemma (\ref{lemma2.3}) we have that the ideal $I_V$ is generated by the single element $\delta(x,x)=\det * [\gamma^1(x_ix_j)]$, where $x=x_1, \ldots ,x_n$ is an $A$-module basis of $E=V$. By Lemma (\ref{lemma2.2}) we have the identity
$ \det[\gamma^1(x_ix_j)] =\det [\gamma^1(x_ix_j)* \gamma^{n-1}(1)]$ in
$\Gamma^n_AF$. As $\sigma_E$ is an algebra homomorphism we have
$$
\sigma_E \det [\gamma^1(x_ix_j)* \gamma^{n-1}(1)]=\det
[\sigma_E(\gamma^1(x_ix_j)* \gamma^{n-1}(1))].$$
By Proposition (\ref{prop1.11}) we have $\sigma_E (\gamma^1(e_ie_j)* \gamma^{n-1}(1))=\operatorname{Trace}(e\mapsto
ex_ix_j)$. Thus we have a matrix with entries $\operatorname{Trace}(e\mapsto
ex_ix_j)$, and the determinant is then the discriminant. 
\end{proof}

\section{Connection with symmetric tensors}

\subsection{A norm vector} Let $F$ be an $A$-algebra, and let $T^n_AF=F\otimes_A \cdots \otimes_AF$ be the tensor product with $n$-copies of $F$. We fix the positive integer $n$, and for any element $x\in F$ we use the following notation
\eqbeg
x_{[j]}=1\otimes \cdots \otimes x \otimes \cdots \otimes 1, \label{411}
\eqend
where the $x$ occurs at the $j$'th component of $T^n_AF$. The group $\Sgot _n$ of $n$-letters acts on $T^n_AF$ by permuting the factors. For any $n$-elements $x=x_1, \ldots , x_n$ in $F$ we define the norm vector
$$\nu (x)=\nu (x_1, \ldots, x_n) =\sum_{\sigma \in \Sgot_n}(-1)^{|\sigma|}x_{\sigma(1)}\otimes \cdots \otimes x_{\sigma (n)} \in T^n_AF,$$
where the summation runs over the elements $\sigma $ in the group $\Sgot_n$. The norm vector $\nu (x)$ is not a symmetric tensor since clearly we have that $\nu (x)=\det (x_{i,[j]})$, where $(x_{i,[j]})$ is the square matrix with coefficients $x_{i,[j]}$ on row $i$, and column $j$, with $1\leq i,j\leq n$. It follows  that
$
\nu(x_{\rho (1)}, \ldots , x_{\rho(n)})=(-1)^{|\rho|}\nu (x_1, \ldots, x_n)$ for  any $\rho \in \Sgot_n $, and that
\eqbeg
\nu (x)=0 \quad \text {if }x_i=x_j, \text{ with }i\neq j.\label{412}
\eqend

\subsection{} Let  $\TS^n_AF$ denote the invariant ring
of $\T^n_AF$ by the natural action of the symmetric group $\Sgot_n$ in
$n$-letters that permutes the factors. We have the map $F \ra \T^n_AF$
sending $x \mapsto x\otimes \cdots \otimes x$, and it is clear that
the map factors through the invariant ring $\TS^n_AF$. The map $F \ra
\TS^n_AF$ determines a norm of degree $n$, as one readily verifies,
hence there exist an $A$-algebra homomorphism 
\eqbeg  \alpha_n \colon
\Gamma^n_AF \ra \TS^n_AF\label{def_alpha} \eqend
 such that $\alpha_n (\gamma^n(x))=x\otimes \cdots \otimes x$, for all
$x
\in F$. 


\subsection{The shuffle product} When $F$ is an $A$-algebra that is flat as an $A$-module, or if $n!$ is invertible in $A$, then
the $A$-algebra homomorphism $\alpha_n$ (\ref{def_alpha}) is an
isomorphism (\cite[IV, \S 5. Proposition IV.5]{roby_lois_pol_mod}, \cite[Exercise 8(a), AIV. p.89]{bourbaki_polyetfractions}). In those cases we can identify $\Gamma_AF$ as the graded sub-module
$$ \Gamma_AF =\bigoplus_{n\geq 0}TS^n_AF \subseteq \bigoplus_{n\geq 0}
\T^n_AF=\T_AF.$$ The external product structure on $\Gamma_AF$ is then
identified with the shuffle product on the full tensor algebra
$\T_AF$. The shuffle product of an $n$-tensor $x\otimes \cdots \otimes
x$ and an $m$-tensor $y\otimes \cdots \otimes y$ is the $m+n$-tensor
given as the sum of all possible different shuffles of the $n$ copies
of $x$ and $m$ copies of $y$ (\cite[Exercise 8
(b), AIV. p.89]{bourbaki_polyetfractions}).

\begin{prop}\label{prop4.4} Let $F$ be an $A$-algebra, and let
$x=x_1, \ldots, x_n$ and $y=y_1, \ldots ,y_n$ be $2n$ elements of $F$. The $A$-algebra homomorphism $\alpha_n \colon \Gamma^n_AF \ra \TS^n_AF$ (\ref{def_alpha}) is such that
$$ \alpha_n(\delta(x,y))=\nu(x)\nu(y).$$
\end{prop}

\begin{proof} The homomorphism $\alpha_n \colon \Gamma^n_AF \ra\TS^n_AF$ takes the
external product $*$ of $\Gamma_AF$ to the shuffle product  of
$\TS_AF$. We then have
\begin{align}
\alpha_n(\delta(x,y))&=\sum_{\sigma \in
  \Sgot_n}(-1)^{|\sigma|}\alpha_n(\gamma^1(x_1y_{\sigma(1)})*\dots
*\gamma^1(x_ny_{\sigma(n)})) \nonumber \\
&=\sum_{\sigma \in
  \Sgot_n}(-1)^{|\sigma|}\sum_{\tau \in
  \Sgot_n}(x_{\tau(1)}y_{\sigma(\tau(1))})\otimes \dots \otimes (x_{\tau(n)}y_{\sigma(\tau(n))}) \nonumber
\end{align}
By formally manipulating the expression above we obtain 
\begin{align}
&=\sum_{\tau \in
  \Sgot_n}(-1)^{|\tau|}\sum_{\sigma \in
  \Sgot_n}(-1)^{|\sigma||\tau|}(x_{\tau(1)}y_{\sigma(\tau(1))})\otimes \dots \otimes (x_{\tau(n)}y_{\sigma(\tau(n))}) \nonumber \\
&=\sum_{\tau \in
  \Sgot_n}(-1)^{|\tau|}\sum_{\tau' \in
  \Sgot_n}(-1)^{|\tau'|}(x_{\tau(1)}y_{\tau'(1)})\otimes
\dots \otimes (x_{\tau(n)}y_{\tau'(n)})\nonumber \\
&=\Big(\sum_{\tau \in \Sgot_n}(-1)^{|\tau|}x_{\tau (1)}\otimes \dots \otimes
x_{\tau (n)}\Big)\Big(\sum_{\tau' \in \Sgot_n}(-1)^{|\tau'|}y_{\tau' (1)}\otimes \dots \otimes
y_{\tau' (n)}\Big)\nonumber \\
&=\nu(x)\nu(y), \nonumber
\end{align} 
and we arrive at $\alpha_n(\delta(x,y))=\nu(x)\nu(y)$ as claimed. 
\end{proof}

\begin{cor}\label{support} Let $\tilde{\alpha} \colon \Gamma^n_AF \ra \T^n_AF$ denote the composition of the map $\alpha_n$  and the inclusion $\TS^n_AF \subseteq \T^n_AF$. Let $I\subseteq \T^n_AF$ denote the extension of the ideal of norms $I_F$ by $\tilde{\alpha}$, and let $J\subseteq \T^n_AF$ denote the ideal of the diagonals. Then we have $\sqrt{I}=\sqrt{J}$.
\end{cor}

\begin{proof} Let $\varphi \colon \T^n_AF \ra L$ be a morphism with $L$ a field, and let $\varphi_i \colon F \ra L$ be the composition of $\varphi$ and the $i$'th co-projection $F \ra \T^n_AF$, where $i=1, \ldots, n$. If $\varphi$ corresponds to a point in the open complement of the diagonals then all the maps $\varphi_i$ are different. That is, no $\p_i =\ker(\varphi_i)$ is contained in another $\p_j$. Furthermore, since the kernels also are prime ideals there exists, for each $i$, an element $x_i$ not in $\p_i$, but where $x_i \in \p_j$ when $j\neq i$. We then have that $\varphi_j(x_i)=0$ for $j\neq 0$, and that $\varphi_i(x_i)\neq 0$. Hence there are elements $x_1, \ldots, x_n$ in $F$ such that $\det (\varphi_j(x_i)) \neq 0$. Then also the image of $\nu(x_1, \ldots, x_n)$ is non-zero in $L$, and we have that the point $\varphi$ is in the open complement of the scheme defined by $I\subseteq \T^n_AF$.

Conversely, if $\varphi$ corresponds to a point on the diagonals then at least two of the maps $\varphi_i$ are equal. Consequently, for any elements $x_1, \ldots, x_n$ in $F$ we have that $\varphi(\nu (x_1, \ldots, x_n))=0$. It follows that $I\subseteq \ker \varphi $, proving the claim.
\end{proof}

\section{Grothendieck-Deligne norm map}

In this section we recall the Grothendieck-Deligne norm map following Deligne (\cite{sga4_deligne_coh_supp_prop}), and we discuss briefly the related Hilbert-Chow morphism. Furthermore we define the notion of sufficiently big sub-modules.

\subsection{The Hilbert functor of $n$-points} We fix an $A$-algebra $F$, and a positive integer $n$. We let $\Hilb^n_F$ denote the covariant functor from the category of $A$-algebras to sets, that sends an $A$-algebra $B$ to the set
\begin{equation*}
\begin{aligned}
{\Hilb}_F^n(B)= \{ & \text{ideals in $F\otimes_A B$ such that the quotient $E$ is} \\
 & \text{ locally free of rank $n$ as a $B$-module} \}.
\end{aligned}
\end{equation*}

\subsection{The Grothendieck-Deligne norm} If $E$ is an $B$-valued point of ${\Hilb}^n_F$ we have the sequence
\eqbeg
F \ra F\otimes_A B \ra E,
\label{321}
\eqend
from where we obtain the $A$-algebra homomorphisms $\Gamma^n_AF \ra \Gamma^n_BE$ that sends $\gamma^n(x)$ to $\gamma^n(\bar x \otimes 1)$, where $\bar x \otimes 1$ is the residue class of $x\otimes 1$ in $E$. Furthermore, when we compose the homomorphism $\Gamma^n_AF \ra \Gamma^n_BE$ with the canonical homomorphism $\sigma_E \colon \Gamma^n_BE \ra B$ we obtain an assignment that is functorial in $B$;  that is we have a morphism of functors
\eqbeg \label{transformationn_F}
\n_F \colon {\Hilb}^n_F \ra \operatorname{Hom}_{A\text{-alg}}(\Gamma^n_AF, -).
\eqend
The natural transformation $\n_F$ we call the Grothendieck-Deligne norm map.

\begin{rem} The Hilbert functor ${\Hilb}^n_F$ can in a natural way be viewed as a contra-variant functor from the category of schemes (over $\Spec(A)$) to sets. In that case the functor ${\Hilb}^n_F$ is representable by a scheme (see e.g.  \cite{GLS_elementary}). If  $X=\Spec (F) \ra S=\Spec(A)$ we write
$\n_X\colon \operatorname{Hilb}^n_{X/S} \ra \Spec (\Gamma^n_AF) $
for the morphism that corresponds to the natural transformation (\ref{transformationn_F}). 
\end{rem}

\subsection{The geometric action}\label{sec1.10} Let $A=K$ be an algebraically closed field, and let $E$ be a finitely generated Artinian $K$-algebra. As $E$ is Artinian it is a product of local rings $E=\prod_{i=1}^pE_i$, and we let $\rho_i \colon E \ra K$ denote the residue class map that factors via $E_i$. Let $m_i=\text{dim}_K(E_i)$, and let $n=\text{dim}_K(E)=m_1+\dots +m_p$. Iversen (\cite[Proposition 4.7]{iversen_lineardeterminants}) shows that the canonical homomorphism $\sigma_E \colon \Gamma^n_KE=\TS^n_KE \ra K$ factors via the homomorphism
$ \rho  \colon T^n_KE \ra K$, where
$$\rho =(\rho_1, \ldots , \rho_1,  \ldots , \rho_p, \ldots ,\rho_p),$$
and where each factor $\rho_i$ is repeated $m_i$-times.

\subsection{Hilbert-Chow morphism} 
Assume that the base ring $A=K$ is a field, and let $X=\Spec(F)$. Then we can
identify $\Spec(\Gamma^n_KF)$ with the symmetric quotient
$\operatorname{Sym}^n(X):=\Spec (\TS^n_KF)$.  Furthermore we have that the
$\Spec(K)$-valued points of $\operatorname{Hilb}^n_X$ correspond to closed
zero-dimensional subschemes $Z\subseteq X$ of length $n$. When $K$ is
algebraically closed we have by (\ref{sec1.10}) that the Grothendieck-Deligne
norm map sends an $K$-valued point $Z\subseteq X$ to the ``associated''
zero-dimensional cycle
$$ \n_X(Z)=\sum_{P\in |Z|}\text{dim}_K(\CalO_{Z,P})[P],$$
where the summation runs over the points in the support of $Z$. Hence we see
that the norm morphism $\n_X$ has the same effect on geometric points as the
Hilbert-Chow morphism. The Hilbert-Chow morphism that appears in
\cite{fogarty_algfam} and \cite{FGAexplained} requires that the Hilbert scheme
is reduced, whereas the Hilbert-Chow morphism that appears in
\cite{kollar_rat_curves} requires that the Hilbert scheme is (semi-) normal. As
the morphism $\n_X$ does not require any hypothesis on the source we have chosen
to refer to that morphism with a different name; the Grothendieck-Deligne norm
map.

\begin{lem}\label{lemma1.6} Let $A=K$ be a field of characteristic zero, and let $F=K[T]$ be the polynomial ring in a finite set of variables $T_1, \ldots ,T_r$. For $n>0$ the $K$-algebra $\Gamma^n_KF$ is generated by
$$ \gamma^1(m)*\gamma^{n-1}(1),$$
for monomials $m\in K[T]$ of degree $\text{deg}(m)\leq n$.
\end{lem}

\begin{proof} The identification $\alpha_n \colon \Gamma^n_KK[T] \ra \TS^n_KK[T]$ identifies, for any $m\in K[T]$, the element $\gamma^1(m)*\gamma^{n-1}(1)$ with the shuffled product of $\alpha_1(m)=m$ and $\alpha_{n-1}(1)=1\otimes \cdots \otimes 1$. That is
$$
\alpha_n(\gamma^1(m)*\gamma^{n-1}(1)) =m\otimes 1\cdots \otimes 1+ \cdots +1\otimes \cdots 1\otimes m=P(m).$$
By a well-known result of Weyl (\cite[II 3]{weyl_invariants}) the invariant ring
$\TS^n_KF$ is generated by the power sums $P(m)$ of monomials $m\in
K[T]$ of degree less or equal to $n$.
\end{proof}


\begin{defn}[Sufficiently big modules] Let us fix an $A$-algebra $F$. An $A$-submodule $V\subseteq F$ is $n$-sufficiently big if the composite $B$-module homomorphism
$$  V\otimes_A B \ra F\otimes_AB \ra E  $$
is surjective for all $A$-algebras $B$, and all $B$-valued points $E$ of the
Hilbert functor ${\Hilb}^n_F$.
\end{defn}

\begin{rem} Sufficiently big submodules always exist as we can take $V=F$.
\end{rem}

\begin{rem} If $V$ is sufficiently big then we clearly have a morphism of functors
$$ {\Hilb}^n_F \ra \mathrm{Grass}^n_V $$
from the Hilbert functor of rank $n$-families, to the Grassmannian of locally free rank $n$-quotients of $V$.
\end{rem}


\begin{thm}\label{thm3.6} Let $F$ be an $A$-algebra, $n$ an positive integer and let $V\subseteq F$ be an $n$-sufficiently big submodule. Then we have for any $A$-algebra $B$, and any $B$-valued point $E$ of ${\Hilb}^n_F$ that the extension of $I_V$, the ideal of norms associated to $V$, by the Grothendieck-Deligne norm map $\n_F \colon \Gamma^n_AF \ra B$ is the discriminant ideal of $E$ over $B$. That is 
$$ \n_F(I_V)B=D_{E/B}\subseteq B.$$
\end{thm}

\begin{proof} We first make a reduction to the situation when $B$ is a local ring. For any $A$-algebra homomorphism  $B\ra B'$ we let $E'=E\otimes_BB'$, and we have the commutative diagram
$$\xymatrix{
\Gamma^n_AF  \ar[r] & \Gamma^n_BE \ar[d] \ar[r]^{\sigma_E}& B\ar[d] \\
 & \Gamma^n_{B'}E' \ar[r]^{\sigma_{E'}} & B'.} $$
Thus by letting $B'$ be the local ring of a prime ideal in $B$, and using the fact that discriminant ideals are compatible with base change, we may assume that $B$ is a local ring.

Let $B$ be a local ring with residue field $K$. As the composition map
$$\xymatrix{
V\otimes_AK \ar[r] & F\otimes_AK \ar[r] & E\otimes_AK}
$$
is surjection of $K$-vector spaces we can find elements $x_1, \ldots, x_n$ in $V$ such that the residue classes of $x_1\otimes \id_K, \ldots , x_n\otimes\id_K$ in $E\otimes_AK$ form a $K$-vector space basis. It then follows from Nakayama's Lemma that the residue classes of $x_1\otimes \id_B,\ldots ,x_n\otimes \id_B$ form a $B$-module basis of $E\otimes_AB=E$. Then we have that the residue class $\bar y \otimes \id_B$ for any element $y\in F$, can be written as a $B$-linear combination of $\bar x_1\otimes \id_B, \ldots ,\bar x_n\otimes \id_B$. In particular we have by Lemma (\ref{lemma2.3}) that $I_{E}\subseteq \Gamma^n_BE$ the ideal of norms associated to $V=E$ is generated by the single element 
\eqbeg \delta (x,x)=\det *[\gamma^{1}(\bar{x_ix_j}\otimes \id_B)]_{1\leq i,j\leq n} \in \Gamma^n_BE.\label{361}
\eqend
By Lemma (\ref{lemma2.5}) we have that the extension of $I_V\subseteq \Gamma^n_AF$ by the composite map $\Gamma^n_AF \ra \Gamma^n_BE$ is exactly $I_E$. Thus we obtain that the extension of $I_V$ by the Grothendieck-Deligne norm map $\n_F$ is the ideal generated by the image of the element (\ref{361}) in $B$. By the definition of $\n_F$ we apply the canonical homomorphism $\sigma_E \colon \Gamma^n_BE \ra B$ to obtain the image of (\ref{361}) in $B$. The result now follows from Proposition (\ref{prop2.7}). 
\end{proof}

\section{Families of distinct points}\label{Sec5}

\subsection{The canonical morphism} The map $ F \ra F\otimes_A \Gamma^{n-1}_AF$
sending $z$ to $z\otimes \gamma^{n-1}(z)$ determines a norm of degree
$n$. Consequently there is a unique $A$-algebra homomorphism $\Gamma^n_AF \ra F
\otimes_A \Gamma^{n-1}_AF $ that takes $\gamma^n(z)$ to $z\otimes
\gamma^{n-1}(z)$. Let
\eqbeg
\pi_n \colon\Spec (F)\times_{\Spec (A)}\Spec (\Gamma^{n-1}_AF) \ra \Spec (\Gamma^n_AF) \label{451}
\eqend
denote the corresponding morphism of schemes. Furthermore, we let $\Delta \subseteq \Spec( \Gamma^n_AF)$ denote the closed subscheme corresponding to the ideal of norms associated to $F$.

\begin{prop}\label{newprop4.8} Let $U=\Spec(\Gamma^n_AF)\setminus \Delta$ denote the open set where the ideal sheaf of norms vanishes. Then the induced morphism
$$ \pi_{n\mid} \colon \pi_n^{-1}(U) \ra U$$
is \'etale of rank $n$.
\end{prop}

\begin{proof} 
Let $U_n\subseteq \Spec(\T^n_AF)$ denote the open complement of the
diagonals. The group of $n$-letters, $\Sgot_n$, acts freely on $U_n$ and the
quotient map $U_n \ra U_n/\Sgot_n$ is \'etale of rank $n!=|\Sgot_n|$. The
morphism $\Spec(\alpha_n) \colon \Spec (\TS^n_AF) \ra \Spec (\Gamma^n_AF)$ is an
isomorphism when restricted to $U_n/\Sgot_n$ (see
e.g. \cite[Prop. 4.2.6]{rydh_fams}). It then follows from Corollary
(\ref{support}) that the map $\tilde{\alpha}_n : \Spec (\T^n_AF) \ra \Spec
(\Gamma^n_AF)$ is an $\Sgot_n$-torsor over $\Spec (\Gamma^n_AF)\setminus
\Delta$. Furthermore, after a faithfully flat base change $A\ra A'$ we can
assume that $\Gamma^n_A(F)\otimes_AA'=\Gamma^n_{A'}(F\otimes_AA')$ is generated
by elements of the form $\gamma^n(z)$ (\cite[Lemma 2.3.1]{ferrand_norme}). Then
clearly the diagram
$$ \xymatrix{
\Spec (\T^n_AF) \ar[rr]^{\tilde{\alpha}_n} \ar[dr]_{1\times
\tilde{\alpha}_{n-1}} && \Spec (\Gamma^n_AF) \\ &\Spec (F) \times \Spec
(\Gamma^{n-1}_AF) \ar[ur]^{\pi_n} }
$$
is commutative.  Over the complement of $\Delta \subseteq \Spec(\Gamma^n_AF)$ we
have that $\tilde{\alpha}_n$ is \'etale of rank $n!$, and $1\times
\tilde{\alpha}_{n-1}$ is \'etale of rank $(n-1)!$. And consequently $\pi_n$ is
\'etale of rank $n$.
\end{proof}

\subsection{Notation} We have the ordered sequence $x=x_1, \ldots, x_n$ of
elements in $F$ fixed. Let $U_A(x)$ be $\Gamma^n_AF$ localized at the element
$\delta(x,x)$, and consider the induced map
\eqbeg
 U_A(x) \ra  (F\otimes_A\Gamma_A^{n-1}F)\otimes_{\Gamma^n_AF}U_A(x)=M_A(x) \label{471}
\eqend
obtained by localization of (\ref{451}).

\begin{lem}\label{lemma4.10} The images of the elements $x=x_1, \ldots ,x_n$ by the map $F \ra F\otimes_A U_A(x) \ra M_A(x)$ form an $U_A(x)$-module basis for $M_A(x)$.
\end{lem}

\begin{proof} By Proposition (\ref{newprop4.8}) we have that $M_A(x)$ is Zariski locally free of rank $n$ over $U_A(x)$. To show that $M_A(x)$ is free it suffices to show that the images of $x_1, \ldots ,x_n$ form a basis locally.  Hence we may assume that $M$ is a free $U$-module, where $U$ is some localization of $U_A(x)$. Let $e=e_1, \ldots ,e_n$ be a basis of $M$, and let $q(z)$ denote the image of $z\in F$ in $M$. There exist scalars $a_{i,j} \in U$ such that 
$q(x_i)=\sum_{j=1}^na_{i,j}e_j$ for $i=1, \ldots, n$. Let $q(x)=q(x_1), \ldots,
q(x_n)$, and let $A=(a_{i,j})$ denote the matrix of the scalars. By Lemma
(\ref{lemma2.3}) we obtain
$$ \delta (q(x),q(x))=\det(A^2) \delta (e,e) \quad \text{in } \quad \Gamma^n_UM.$$
 The element $\delta(x,x)\otimes 1$ in $\Gamma^n_A(F)\otimes_{\Gamma^n_AF}U=\Gamma^n_U(F\otimes_AU)$ is invertible by definition. The natural morphism $F\otimes_AU \ra M$ induces a morphism $\Gamma^n_U(F\otimes_AU) \ra \Gamma^n_UM$ sending $\delta(x,x)\otimes 1$ to the invertible element $\delta(q(x),q(x))$. Then also $\det (A)$ must be invertible, and consequently we have that $q(x_1), \ldots, q(x_n)$ form a basis of $M_A(x)$.
\end{proof}

\begin{defn} The functor $\CalH^{et}_F(x)$ is the covariant functor from the category of $A$-algebras to sets that map an $A$-algebra $B$ to the set of ideals in $F\otimes _A B$ such that corresponding quotients $Q$ satisfy the following
\renewcommand{\labelenumi}{(\theenumi)}
\begin{enumerate}
\item The elements $q(x_1), \ldots ,q(x_n)$ in $Q$ form a $B$-module basis, where $q : F \ra F\otimes_AB \ra Q$ is the composite map.
\item The algebra homomorphism $ B \ra Q$ is \'etale.
\end{enumerate}
\end{defn}

\begin{lem}\label{lemma5.2} Let $B$ be an $A$-algebra, and $Q$ a $B$-valued point of $\CalH^{et}_F(x)$. Then we have the following commutative diagram of algebras
\eqbeg \label{521}
\xymatrix{
\Gamma^n_AF \ar[r] \ar[d]^{\operatorname{can}} & \Gamma^n_BQ \ar[d]^{\sigma_Q} \\
U_A(x)\colon=(\Gamma^n_AF)_{\delta(x,x)} \ar[r] & B.}
\eqend
\end{lem}

\begin{proof} The composite morphism $F \ra F\otimes _AB \ra Q$ induces a morphism of $A$-algebras $\Gamma^n_AF \ra \Gamma^n_BQ$ that sends the element $\delta(x,x)$ to $\delta (q(x),q(x))$, where $q(x)=q(x_1), \ldots ,q(x_n)$ in $Q$. By assumption the elements $q(x)$ form a basis of $Q$ and that  $Q$ is \'etale. Then, by Proposition (\ref{prop2.7}) we that the image of $\delta (q(x),q(x))$ by the canonical map $\sigma_Q \colon \Gamma^n_BQ \ra Q$ is a unit, and the commutativity of the diagram (\ref{521}) follows.  
\end{proof}

\subsection{Universal coefficients}  For each pair of indices $i\leq j$ we look at the product $x_ix_j$ in $F$, and for each $k=1, \ldots, n$ we consider the sequence
\eqbeg x^{i,j}_k =x_1, \ldots , x_{k-1}, x_{i}x_j, x_{k+1}, \ldots ,x_n\label{541}
\eqend
where the $k$'th element is replaced with the product $x_ix_j$. We now define the universal coefficient
\eqbeg
 \alpha^{i,j}_k =\frac{\delta(x,x^{i,j}_k)}{\delta (x,x)} \quad \text{in}\quad U_A(x)=(\Gamma^n_AF)_{\delta(x,x)}.\label{542}
\eqend

\begin{prop}\label{Repl.prop5.5} Let $Q$ be a $B$-valued point of $\CalH^{et}_F(x)$, and let $q\colon F \ra F\otimes_AB \ra Q$ denote the composite map. Let $b^{i,j}_k$ be the unique elements in $B$ such that
$$ q(x_ix_j) =\sum_{k=1}^n b^{i,j}_k q(x_k)$$
in $Q$. Then $b^{i,j}_k$ is the specialization of the element $\alpha^{i,j}_k$ under the natural map $U_A(x) \ra B$ of Lemma (\ref{lemma5.2}), for each $i,j,k =1, \ldots, n$. In particular we have that $M_A(x)\otimes_{U_A(x)}B=Q$ as quotients of $F\otimes_AB$.
\end{prop}

\begin{proof} Having the triplet $i,j,k$ fixed, we let  $x^{i,j}_k$ denote the sequence (\ref{541}) of elements in $F$. Consider the element $\delta(q(x),q(x^{i,j}_k))$ in $\Gamma^n_BQ$. We replace the element $q(x_ix_j)$ in $Q$ with $\sum b_k^{i,j}q(x_k)$, and obtain 
$$ 
\delta(q(x),q(x^{i,j}_k)) = b_k^{i,j}\delta(q(x),q(x)) \quad \in \Gamma^n_BQ.
$$
The element $\delta(q(x),q(x))$ is the image of $\delta(x,x)$ by the induced map $\Gamma^n_AF \ra \Gamma^n_BQ$. It follows from the commutative diagram (\ref{521}) that $b_k^{i,j}$ in $B$ is the image of $\alpha^{i,j}_k$.
\end{proof}

\begin{cor}\label{reprH^et(x)} The pair $(U_A(x),M_A(x))$ represents $\CalH^{et}_F(x)$.
\end{cor}

\begin{proof} It follows from Proposition (\ref{newprop4.8}) and Lemma (\ref{lemma4.10}) that $M:=M_A(x)$ is a $U:=U_A(x)$-valued point of $\CalH^{et}_F(x)$. If $Q$ is any $B$-valued point of  $\CalH^{et}_F(x)$ have by Proposition (\ref{Repl.prop5.5}) one morphism $U \ra B$ with the desired property, and we need to establish uniqueness of that map. Therefore, let $\varphi _i : U \ra B$ ($i=1,2$), be two $A$-algebra homomorphisms such that both extensions $M\otimes_{U}B$ equal $Q$ as quotients of $F\otimes_AB$. We then have that the natural map 
$$\Gamma^n_{U}M \ra\Gamma^n_{U}(M)\otimes_{U}B=\Gamma^n_BQ$$
is independent of the maps $\varphi_i : U \ra B$. And in particular the canonical section $\sigma_Q =\sigma_M\otimes 1 : \Gamma^n_BQ \ra B$ is independent of the maps $\varphi_i, (i=1,2)$. For any element $u\in U$ we have that $\sigma_M (u\gamma^n(1))=u$, and then also that $\sigma_Q(u\gamma^n(1)\otimes 1_B)=\varphi_i(u)$. Thus $\varphi_1=\varphi_2$, and we have proven uniqueness.
\end{proof}

\subsection{\'Etale families} We let $\CalH^{et,n}_{F}$ denote the functor of \'{e}tale families of the Hilbert functor ${\Hilb}^n_{F}$ of $n$-points on $F$. That is, we consider the co-variant functor from $A$-algebras to sets whose $B$-valued points are
$$ \CalH^{et,n}_{F}(B) =\{ I \in {\Hilb}^n_F(B) \mid B \ra F\otimes_AB/I \text{ is \'etale}\}.$$
It is clear that $\CalH^{et,n}_{F}$ is an open subfunctor of ${\Hilb}^n_{F}$ and we will end this section by describing the corresponding open subscheme of the Hilbert scheme.

\begin{prop}\label{prop5.9} Let $F$ be an $A$-algebra. Let $\Delta\subseteq \Spec(\Gamma^n_AF)$ be the closed subscheme defined by the ideal of norms $I_F$, and and let $U=\Spec(\Gamma^n_AF)\setminus \Delta$ denote its open complement. The family
$ \pi_{n\mid} \colon \pi_n^{-1}(U) \ra U$
given in (\ref{451}) represents $\CalH^{et,n}_{F}$.
\end{prop}

\begin{proof} Clearly the functors $\CalH_F^{et}(x)$, for different choices of elements $x=x_1, \ldots, x_n$ in $F$, give an open cover of $\CalH_F^{et,n}$. By Corollary (\ref{reprH^et(x)}) the restriction of the family $\pi_{n\mid} \colon \pi_n^{-1}(U) \ra U$ to the open subscheme $\Spec (U_A(x))\subseteq U$ represents $\CalH_F^{et}(x)$. Thus we need only to check that union of the schemes $\Spec(U_A(x))$, for different $x=x_1, \ldots, x_n$, is $U$, which however is a consequence of Lemma (\ref{product}).
\end{proof}

\section{Closure of the locus of distinct points}

We will continue with the notation from the preceding sections. In this section
we will construct universal families, not for the locus of distinct points as
in Section \ref{Sec5}, but for its closure.

\subsection{Notation} Let $F$ be an $A$-algebra, and let $R=\oplus_{m\geq 0}I_F^m$ denote the graded ring where
$I_F\subseteq \Gamma^n_AF$ is the ideal of norms associated to $V=F$. We let
$x=x_1, \ldots ,x_n$ be $n$-elements in $F$, and we denote by
$R(x)=R_{(\delta(x,x))}$ the degree zero part of the localization of $R$ at
$\delta(x,x)\in I_F$. Finally we let $\CalE$ denote the free $R(x)$-module of
rank $n$. We will write
\eqbeg
 \CalE = \oplus_{i=1}^nR(x)[x_i],\label{611}
\eqend
where $[x_i]$ is our notation for a basis element pointing out the $i$'th component of the direct sum $\CalE$. As $\Gamma^n_AF$ is an $A$-algebra we have that $\CalE$ is an $A$-module. We define the $A$-module homomorphism
\eqbeg
 [ \ ] \colon F \ra \CalE\label{612}
\eqend
in the following way. For any $y\in F$, and any $i=1, \ldots ,n$, we let
\eqbeg x^i_y=x_1, \ldots ,x_{i-1}, y, x_{i+1}, \ldots ,x_n\label{613}
\eqend
denote the $n$-elements in $F$ where the $i$'th element $x_i$ is replaced with $y$. Then we define the value of the map (6.1.3) on the element $y\in F$ as
\eqbeg
[y]=\sum_{i=1}^n\frac{\delta(x,x_y^i)}{\delta(x,x)}[x_i] \quad \text{ in } \quad \CalE.\label{614}
\eqend
Note that when $y=x_i$ the notation of (\ref{611}) is consistent with the notation of (\ref{614}). 
As determinants are linear in its columns (and rows) it follows that the map $[ \ ] \colon F \ra \CalE$ defined above is an $A$-module homomorphism. Furthermore it follows that we have an induced $R(x)$-module homomorphism
\eqbeg
    [\ ] \otimes \id \colon F\otimes_A R(x) \ra \CalE\otimes_A R(x)\cong\CalE,\label{615}
\eqend
which is surjective.

\subsection{Universal multiplication}\label{subsec6.2} With the notation as above we define now the $R(x)$-bilinear map $\CalE\times \CalE \ra \CalE$ by defining its action on the basis as
\eqbeg
[x_i][x_j] :=[x_ix_j] \quad \text{ for } \quad i,j \in \{1, \ldots, n\}. \label{621}
\eqend
We will show that the above defined bilinear map defines as multiplication structure on $\CalE$ - that is giving $\CalE$ a structure of a commutative $R(x)$-algebra. We first observe the following simple but important fact. Consider $\CalE$ as a sheaf on $\Spec (R(x))$, and let $U\subset \Spec (R(x))$ be a subscheme of $\Spec(R(x))$. Assume furthermore that the bilinear map (\ref{621}) restricted to $\CalE_U$ gives a ring structure on $\CalE_U$. That is the product (\ref{621}) is associative, has an multiplicative identity and is distributive, then we also have a ring structure on $\CalE_{\bar U}$, where $\bar U$ is the scheme theoretic closure of $U\subseteq \Spec (R(x))$. We will apply this observation to a scheme theoretic dense open subset $U\subseteq \Spec (R(x))$.

\begin{prop}\label{prop6.3} Let $F$ be an $A$-algebra. We have that (\ref{614}) defines an algebra structure on $\CalE$ and that the map (\ref{615}) is a surjective $R(x)$-algebra homomorphism.
\end{prop}

\begin{proof} Let $R=\oplus _{n\geq 0}I^n_F$, where $I_F \subseteq \Gamma^n_AF$ is the ideal of norms. We have that $\Spec (R(x))$ is an affine open subset of $\Proj (R)$, where 
$$\rho \colon \Proj (R) \ra \Spec (\Gamma^n_AF)$$
is the blow-up with center $\Delta  =\Spec (\Gamma^n_A(F/I_F))$. The open complement $\Proj (R) \setminus \rho ^{-1}(\Delta)$ of the effective Cartier divisor $\rho ^{-1}(\Delta)$ is schematically dense. Hence
$$ U:= \Spec(R(x)) \setminus \rho^{-1}(\Delta)\cap \Spec (R(x))$$
is schematically dense in $\Spec (R(x))$. By (\ref{subsec6.2}) it suffices to show the statements over $U$. However we have that $U=\Spec (U_A(x))$ as defined in (\ref{521}), and that the restriction of $\CalE_{\mid U}$ coincides with the family $\Spec (M_A(x))$. In other words, we have that restriction of the multiplication map (\ref{615}) to the open $U$ coincides with the universal multiplication map of Proposition (\ref{prop5.9}).
 
\end{proof}

\begin{cor}\label{cor6.4} We have that $\CalE(x)$ is an $R(x)$-valued point of the Hilbert functor ${\Hilb}^n_{F}$.
\end{cor}
\begin{proof} The proposition gives that $\Spec (\CalE)$ is a closed subscheme of $\Spec(F\otimes_A R(x))$. By construction the $R(x)$-module $\CalE$ is free of rank $n$. 
\end{proof}

\begin{cor}\label{cor6.5} 
The schemes $\Spec (R(x))$, for different choices of $x=x_1, \ldots
,x_n$ in $F$, form an affine open cover of $\Proj (R)$, and the
families $\Spec (\CalE(x))\ra \Spec (R(x))$ glue together to a $\Proj
(R)$-valued point of the Hilbert functor ${\Hilb}^n_{F}$.
\end{cor}

\begin{proof} 
The first statement follows from Lemma (\ref{product}). To prove the
second assertion it suffices to see that the families glue over a open
dense set. Let $U=\Proj (R) \setminus \rho^{-1}(\Delta)$, where $\rho
\colon \Proj (R) \ra \Spec(\Gamma^n_AF)$ is the blow-up with center
$\Delta$. Then we have that $\Spec (R(x)) \cap U=\Spec (U_A(x))$ for
any $n$-elements $x=x_1, \ldots ,x_n$ in $F$, and the result follows.
\end{proof}

\section{The good component}

\subsection{} 
When $X\ra S$ is an algebraic space we have the Hilbert functor
$\Hilb^n_{X/S}$ of closed subspaces of $X$ that are flat and finite of
rank $n$ over the base. If $U\ra X$ is an \'etale map we define the
subfunctor $\CalH^n_{U\to X}$ of ${\Hilb}^n_{U/S}$ by assigning to any
$S$-scheme $T$ the set
\begin{equation*}
\begin{aligned}
\CalH^n_{U\to X}(T)= \{ & Z\in {\Hilb}_{U/S}^n(T) \text{ such that the
composite map }\\ & Z \subseteq U\times_ST \ra X\times_ST \text{ is a
closed immersion} \}.
\end{aligned}
\end{equation*}

\begin{prop}\label{eqHilb} 
Let $X \ra S$ be a separated quasi-compact algebraic space over an affine scheme
$S$, and let $U\ra X$ be an \'etale representable cover with $U$ an affine
scheme, and let $R=U\times_XU$. Then we have the following
\renewcommand{\labelenumi}{(\theenumi)}
\begin{enumerate}
\item The functor $\CalH^n_{U\to X}$ is representable by a scheme.
\item The natural map $\CalH^n_{U\to X} \ra {\Hilb}^n_{X/S}$ is
representable, \'etale and surjective.
\item The two maps $ \xymatrix@M=1pt{ \CalH^n_{R\to X} \ar@<.5ex>[r]
\ar@<-.5ex>[r] & \CalH^n_{U\to X}}$ form an \'etale equivalence
relation, and the quotient is $\Hilb^n_{X/S}$.
\end{enumerate}
\end{prop}

\begin{proof} 
Since $X\ra S$ is separated the composition $Z\ra U\times_ST \ra
X\times_ST$ will be finite, for any $Z \in {\Hilb}^n_{U/S}(T)$, any
$S$-scheme $T$. It is then clear that $\CalH^n_{U\to X}$ is an open
subfunctor of ${\Hilb}_{U/S}^n$ where the latter is known to be
representable (\cite{GLS_elementary}). This shows the first assertion.
To see that the map $\CalH^n_{U\to X} \ra \Hilb^n_{X/S}$ is
representable we let $T \ra \Hilb^n_{X/S}$ be a morphism, with $T$
some $S$-scheme. Let $Z\subseteq X\times_ST$ denote the corresponding
closed subscheme, and let $Z_U=Z\times_X U$. It is easily verified
that the fiber product $\CalH^n_{U\to X}\times_{\Hilb^n_{X/S}}T$ equals
the set of sections of $Z_U \ra Z$. Thus the fibred product equals the
Weil restriction of scalars $\R_{Z/T}(Z_U)$ of $Z_U$ with respect to
$Z\ra T$. The fiber of $Z_U \ra T$ over any point in $T$, is an affine
scheme, and it follows from \cite[Thm.4]{BLRaynaud_neron} that the
Weil restriction $\R_{Z/T}(Z_U)$ is representable by a scheme. Hence
the map $\CalH^n_{U\to X} \ra \Hilb^n_{X/S}$ is
representable. \'Etaleness of the map follows from
\cite[Prop. 5]{BLRaynaud_neron}, and surjectivity follows as any
$T$-valued point of $\Hilb^n_{X/S}$ \'etale locally lifts to $U$.  The
last assertion follows as it is easy to see that the natural map
$\CalH^n_{R\to X} \ra \CalH^n_{U\to X}\times_{\Hilb^n_{X/S}}
\CalH^n_{U\to X}$ is in fact an isomorphism.
\end{proof}
\begin{cor} 
Let $X \ra S$ be a separated map of algebraic spaces. Then
${\Hilb}^n_{X/S}$ is an algebraic space.
\end{cor}

\begin{proof} 
It suffices to show the statement for affine base $S$. Let $X'
\subseteq X$ be an open immersion. Then as $X \ra S$ is assumed
separated we have a map ${\Hilb}_{X'/S}^n \ra {\Hilb}_{X/S}^n$ which
is a representable open immersion. Furthermore, as
$$ 
\Hilb^n_{X/S} =  \lim_{
\underset{\text{open, q-compact}}{X'\subseteq X}} 
{\Hilb}_{X'/S}^n
$$ 
we may assume that $X\ra S$ is quasi-compact as well. Then the result
follows from the proposition.
\end{proof}

\begin{rem} 
For a quasi-projective scheme $X \ra S$ over a Noetherian base scheme
$S$ it was proven by Grothendieck that the Hilbert functor
${\Hilb}_{X/S}^n$ is representable by a scheme
(\cite{gr_existence}). For a separated algebraic space $X\ra S$
locally of finite type, Artin proved that ${\Hilb}^n_{X/S}$ is an
algebraic space (\cite{artin_alg_form_moduli_I}). The proof of the
general result above showing that ${\Hilb}^n_{X/S}$ is an algebraic
space for any separated algebraic space $X \ra S$ was suggested to us
by the referee.
\end{rem}

\subsection{The good component} 
Let $X\ra S$ be a separated algebraic space, and let $Z \ra
{\Hilb}^n_{X/S}$ be the universal family, which by definition is
finite, flat of rank $n$. The discriminant $D_{Z}\subseteq
{\Hilb}^n_{X/S}$ is a closed subspace with the open complement
$U^{et}_{X/S}$ parameterizing length $n$ \'etale subspaces of $X$. We
define $\good\subseteq {\Hilb}_{X/S}^n$ as the schematic closure of
the open subspace $U^{et}_{X/S}$. We call $\good $ the {\it good} or
{\it principal component}.

\begin{rem} Let $f\colon Z \ra H$ be a morphism of algebraic spaces which is a
finite and flat morphism of rank $n$. Then the set $U\subseteq H$ where $f$ is
\'etale is an open subset being the complement of the discriminant $D_{Z/H}$.
The scheme theoretic closure of $U\subseteq H$ is then the largest closed
subscheme of $H$ over which the discriminant of $f$ is a non-zero-divisor.
\end{rem}

\begin{thm}\label{thm7.3} 
Let $X=\Spec(F) \ra S=\Spec(A)$ be a morphism of affine schemes, and
let $\Delta \subseteq \Spec(\Gamma^n_AF)$ be the closed subscheme
defined by the ideal of norms. Then we have that the good component
$\good$ is isomorphic to the blow-up $\mathrm{Bl}(\Delta)$ of
$\Spec(\Gamma^n_AF)$ along $\Delta$.  The isomorphism
$$
\xymatrix{ \bn_X \colon \good \ar[r]^{\simeq} & \mathrm{Bl}(\Delta), }
$$ 
is induced from restricting the norm map $\n_X \colon
\operatorname{Hilb}^n_{X/S} \ra \Spec(\Gamma^n_AF)$ to the good
component $\good $.
\end{thm}

\begin{proof}  
By Theorem (\ref{thm3.6}) we have that the inverse image
$\n_X^{-1}(\Delta )$ is the discriminant $D_{Z}\subseteq
\mathrm{Hilb}^n_{X/S}$ of the universal family $Z \ra
\mathrm{Hilb}^n_X$. Consequently we have that the local equation of
the closed immersion
$$ 
\good \cap \n_X^{-1}(\Delta )\subseteq \good ,  
$$ 
is not a zero divisor. Therefore, by the universal properties of the
blow-up, we get an induced morphism $\bn_X \colon \good \ra
\mathrm{Bl}(\Delta)$. A morphism we will show is an isomorphism.

We have by Corollary (\ref{cor6.5}) the $\mathrm{Bl}(\Delta )$-valued point
$\CalE$ of the Hilbert functor ${\Hilb}^n_{F}$. From the defining properties of
the Hilbert scheme we then have a morphism $f_\CalE \colon \mathrm{Bl}(\Delta)
\ra \mathrm{Hilb}^n_{X/S}$ such that the pull-back of the universal family is
$\CalE$.  When restricting $\CalE$ to the open set
$U=\Spec(\Gamma^n_AF)\setminus \Delta $ we have an \'etale family -- by
construction of $\CalE$. Hence the image $f_\CalE(U)$ is contained in
$U^{et}_{X/S}$. It follows that the schematically closure
$\overline{U^{et}_{X/S}}=\good $ contains the image of the closure of
$\overline{U}=\mathrm{Bl}(\Delta)$. Consequently we have a morphism $f_\CalE
\colon \mathrm{Bl}(\Delta) \ra \good $, a morphism we claim is the inverse to
the map $\bn_X \colon\good \ra \mathrm{Bl}(\Delta)$.

By Proposition (\ref{prop5.9}) we have that the restriction of $f_\CalE$ to $U$
is the inverse of the restriction of $\bn_X$ to $U^{et}_{X/S}$.  As both $U$ in
$\mathrm{Bl}(\Delta)$ and $U^{et}_{X/S}$ in $\good$ are open complements of
effective Cartier divisors it follows that $f_\CalE$ is the inverse of $\bn_X$.
\end{proof}

\subsection{}  

For a separated map of algebraic spaces $X \ra S$ there exists an algebraic
space $\Gamma^n_{X/S}$ that naturally globalize the affine situation with $\Spec
(\Gamma^n_{A}F)$ (\cite{rydh_fams}). For the convenience of the reader we will
give a description of this space for $X$ quasi-compact over an affine base. Not
only is the quasi-compact case technically easier to handle, but it turns out
to be sufficient in order to generalize Theorem (\ref{thm7.3}) for separated
algebraic spaces $X \ra S$.

\subsection{Pro-equivalence} 

We will say that two sequences (indexed by the non-negative integers) of ideals
$\{I_m\}$ and $\{J_m\}$ in a ring $B$ are {\em pro-equivalent} if there exists
an integer $p$ such $I_{m+p}\subseteq J_m$, and $J_{m+p}\subseteq I_m$, for all
$m\geq 0$.

\begin{lem}\label{invring} 
Let $G$ be a finite group acting on a Noetherian ring $B$ and let $\A \subseteq
B$ be an invariant ideal. Assume furthermore that the invariant ring $B^G$ is
Noetherian, and that $B$ is a finite module over the invariant ring.  Then
$\{(\A^G)^m\}$ is pro-equivalent with $\{(\A^m)^G\}$.
\end{lem}

\begin{proof} Clearly $(\A^G)^{m+p}\subseteq (\A^{m})^G$, and consequently it suffices to show that $(\A^{m+p})^{G}\subseteq (\A^G)^m$ for some $p$. An element $x\in B$ is a root of the monic polynomial $m_x(t)=\prod_{g\in G}(t-g.x)$. Since $\A$ is $G$-invariant we have for any $x\in \A$ that  $x^{\mid G\mid}\in \A^G$, and consequently that
$$ \A^{m+\mid G \mid} \subseteq (\A^G)^mB.$$
By assumption $B$ is a finitely generated $B^G$-module, and consequently by the Artin-Rees Lemma (\cite[Cor. 10.10]{atiyahmacdonald}) there exists an integer $k$ such that
$$ (\A^G)^mB\cap B^G =(\A^G)^{m-k}\big( (\A^{G})^kB\cap B^G\big) \subseteq (\A^{G})^{m-k}.$$
Hence $(\A^{m+\mid G \mid +k})^G \subseteq (\A^G)^m$ for all $m\geq 0$.
\end{proof}

\begin{lem}\label{proeq} Let $F$ be an $A$-algebra of finite type, and let $I\subseteq F$ be an ideal of finite type. For each $m>0$ we let $J_m$ denote the kernel of the natural map $\Gamma^n_A(F) \ra \Gamma^n_A(F/I^m)$.  Then $\{J_m\}$ is pro-equivalent with $\{J_1^m\}$.
\end{lem}

\begin{proof} 
We first show the case with a polynomial ring over the integers. Let $X=x_1,
\ldots, x_r$ and $T=t_1, \ldots, t_s$ be variables over ${\bf Z}$, and let
$F={\bf Z}[X,T]$, and $I=(T)$. Let $\A_m$ denote the kernel of $\T^n_AF \ra
\T^n_A(F/I^m)$. It is easily checked that $\{ \A_m \}$ is pro-equivalent with
$\{\A_1^m\}$. The group $\Sgot_n$ acts on $\T^n_AF$, and it follows that $\{
(\A_1^m)^{\Sgot_n} \}$ is pro-equivalent with $\{ \A_m^{\Sgot_n}\}$. By Lemma
(\ref{invring}) we have that $\{ (\A_1^m)^{\Sgot_n} \}$ is pro-equivalent with
$\{ (\A_1^{\Sgot_n})^m \}$. As $F/I^m$ is free, an in particular flat ${\bf
Z}$-module for all $m>0$, we have that $\Gamma^n_A(F/I^m)=\TS^n_A(F/I^m)$. In
particular we get that
$$ \ker ( \Gamma^n_AF \ra \Gamma^n_A(F/I^m)) = (\A_m)^{\Sgot_n},$$
and we have proven the lemma in the special case. Since $\Gamma^n_{\bf Z}{\bf
Z}[X,T] \otimes_{\bf Z} A =\Gamma^n_AA[X,T]$ we have also proven the lemma for
$F=A[X,T]$, and $I=(T)$. In the general case we let $\varphi \colon A[X,T]\ra F$
denote the $A$-algebra homomorphism that sends $X$ to a set of generators of
$F$, and $T$ to a set of generators of the ideal $I\subseteq F$. For each $m>0$
we have induced surjective maps $\varphi_m \colon A[X,T]/(T)^m \ra F/I^m$ and
$\Gamma(\varphi_m) \colon \Gamma^n_AA[X,T]/(T)^m \ra \Gamma^n_AF/I^m$. An
element in $\ker(\Gamma(\varphi_m))$ is of the form (\cite[Prop. IV.8,
p. 284]{roby_lois_pol_mod})
$$\gamma^c(\bar{f})*\gamma^{n-c}(\bar{g})$$
where $\bar{g} \in A[X,T]/(T)^m$ and $\bar{f} \in \ker(\varphi_m)$. Clearly we
can find elements $f$ and $g$ in $A[X,T]$, with $f\in \ker(\varphi)$, that
restricts to $\bar{f}$ and $\bar {g}$ by the canonical map. Thus the induced map
$\ker(\Gamma^n(\varphi)) \ra \ker(\Gamma^n(\varphi_m))$ is surjective for all
$m>0$. It follows that the induced map from
$$\A_m=\ker \big( \Gamma^n_AA[X,T] \ra \Gamma^n_A(A[X,T]/(T)^m) \big)$$
to $J_m =\ker(\Gamma^n_AF \ra \Gamma^n_A(F/I^m))$ is surjective. In particular
$\A_1$ surjects to $J_1$, so $\A^m_1$ surjects to $J_1^m$. The lemma now follows
by lifting elements to $\A_m$ and $\A_1^m$, where the result holds.\end{proof}

\subsection{FPR-sets} 
Let $G$ be a finite group acting on a separated algebraic space
$X$. By a result of Deligne the geometric quotient $X/G$ exists as an
algebraic space. We will make use of that result, but we need also to
recall some terminology: For any group element $\sigma \in G$ we have
the induced map $(\id_X,\sigma) \colon X \ra X\times_SX$. By taking
the intersection of the diagonal and $X$ via the map $(\id_X,\sigma)$
we get a closed subspace $X^{\sigma}\subseteq X$. If $f \colon X \ra
Y$ is a $G$ equivariant map, we have a closed immersion $X^{\sigma}
\subseteq f^{-1}(Y^{\sigma})$. An equivariant map $f \colon X \ra Y$
is \emph{fixed-point-reflecting} (abbreviated \fpr) if we have an
equality of {\em sets} $X^{\sigma} =f^{-1}(Y^{\sigma})$ for all
$\sigma \in G$ (\cite[p.183]{knutson_alg_spaces}). We have an
alternative description of this condition: For every point $x \in X$
the stabiliser group $G_x$ is equal to the stabiliser group $G_y$ of
the image point $y=f(x)$ (in general we only have the inclusion
$G_x\subseteq G_y$). This allows us to say that $f$ is
\emph{fixed-point-reflecting at $x$} (abbreviated \fpr{} at $x$) if
$G_x=G_y$. We then have some general facts about \fpr-sets.
\begin{lem}\label{fpr-lemma}
i) Suppose $f\colon X \to Y$ and $g\colon Y \to Z$ are $G$-morphisms,
$h$ their composite and $x \in X$. Then $h$ is \fpr{} at $x$
precisely when $f$ is \fpr{} at $x$ and $g$ is \fpr{} at $f(x)$.

ii) Suppose that $\{X_\alpha\}$ is an inverse system with affine
transition maps of $G$-spaces. For $x \in X := \ili_{\alpha}X_\alpha$
the set $S_x:=\{\alpha\mid p_\alpha \text{ is \fpr{} at }x\}$, where
$p_\alpha\colon X \to X_\alpha$ is the structure map, is non-empty and
\emph{upwards closed} (i.e., if $\alpha \in S_x$ and $\alpha
'\ge\alpha$ then $\alpha '\in S_x$).

iii) Suppose now that also $\{Y_\alpha\}$ is an inverse system with
affine transition maps of $G$-spaces over the same index set and that
$\{f_\alpha\colon X_\alpha \to Y_\alpha\}$ is a $G$-morphism of
directed systems. Set $Y:=\ili_\alpha Y_\alpha$, $f:=\ili_\alpha
f_\alpha$ and assume that $f$ is \fpr{} at $x \in X$. Then
$\{\alpha\mid f_\alpha \text{ is \fpr{} at }x_\alpha\}$ is non-empty
and upwards closed.
\end{lem}
\begin{proof}
For the first part we always have that $G_x \subseteq G_{f(x)} \subseteq
G_{h(x)}$ so that if $h$ is \fpr{} at $x$, i.e., $G_x=G_{h(x)}$, then $f$ is
\fpr{} at $x$ and $g$ is \fpr{} at $f(x)$ and clearly conversely.

That $S_x$ is upwards directed follows from i). For every $g \notin
G_x$ we have $gx\ne x$ and hence there is an index $\alpha_g$ such
that $gx_{\alpha_g}\ne x_{\alpha_g}$, where for all indices $\beta$,
$x_{\beta}:=p_{\beta}(x)$. Picking an $\alpha$ such that $\alpha \ge
\alpha_g$ for all such $g$ we get that $gx_\alpha\ne x_\alpha$ for all
$g \in G_x$ which implies that $G_{x_\alpha}=G_x$ so that $\alpha \in
S_x$.

Finally, iii) follows from i) and ii) applied to $f(x) \in Y$.
\end{proof}

\begin{deflem}
If the equivariant map $f \colon X \ra Y$ is separated and unramified, then $X^{\sigma}$ is
both open and closed in $f^{-1}(Y^{\sigma})$. Hence if $Y$ is also separated
over some $S$ on which $G$ acts trivially, there is a maximal open
\fpr-subset of $X$, which we call the \emph{\fpr-locus} of $f$.

In the particular case when $U\ra X$ is an unramified separated map and $X$ is
separated over $S$, we will denote the \fpr-locus of the induced $\Sgot_n$-map
$U^n_S \ra X^n_S$ by $\Omega_{U\to X}\subseteq U^n_S$.
\end{deflem}
\begin{proof}
We have a map $f^{-1}(Y^{\sigma}) \to X\times_YX$ given by $x \mapsto (x,\sigma
x)$ and $X^\sigma$ is the inverse image of the diagonal. As $f$ is unramified
and separated, the diagonal is open and closed in $X\times_YX$ and hence so is
$X^\sigma$ in $f^{-1}(Y^{\sigma})$. If $Y$ is also separated, then $f^{-1}(Y^\sigma)$ is
closed in $X$ and hence the complement of $X^\sigma$ in $f^{-1}(Y^\sigma)$ is
closed in $X$ and removing such subsets for all $\sigma$ gives the \fpr-locus.
\end{proof}

\begin{lem}\label{inftyngh} 
Let $F\ra F'$ be an \'etale extension of $A$-algebras.  Let $\varphi
\colon \T^n_AF' \ra L$ be a map to a field $L$ and let $\varphi _i
\colon F' \ra L$ be the co-projections of $\varphi$ (with $i=1, \ldots
, n$). Define the ideals $J=\cap \ker \varphi_i $ in $F'$ and $I=\cap
\ker \varphi_{i| F}$ in $F$. If the point $\varphi$ is in the fixed
point reflecting set $\Omega_{F\to F'}$ of $\Spec(\T^n_AF') \ra \Spec
(\T^n_AF)$, then the induced map
$$ F/I^m \ra F'/J^m$$
is an isomorphism, for all $m>0$.
\end{lem} 

\begin{proof}
We may work \'etale locally around the points given by $\ker\varphi_i$
and $\ker\varphi_{i|F}$ so we may assume that $F$ is a product of
strictly Henselian rings with $\ker\varphi_{i|F}$ as maximal ideals
and similarly for $F'$. That $\varphi \in \Omega_{F\to F'}$ means that
the map $F \to F'$ induces a bijection on maximal ideals, which as $F$
is semi-local strictly Henselian and $F \to F'$ is \'etale means that
$F \to F'$ is an isomorphism. From this the lemma follows immediately.
\end{proof}

\subsection{Notation}\label{fpringamma} 
When $U \ra X$ is an \'etale cover, we let
$$
\Omega_{U\to X}'\subseteq U^n_S/\Sgot_n
$$
be the image of the \fpr-locus $\Omega_{U\to X}\subseteq U^n_S$ by the quotient
map. Set $R=U\times_XU$, then the \fpr-locus $\Omega_{R\to X}$ is identified
with $\Omega_{U\to X}\times_{X^n_S}\Omega_{U\to X}$, so in particular we have
that (\cite[p. 183-184]{knutson_alg_spaces}) \eqbeg \label{eteqrelsym}
\xymatrix{ {{\Omega}_{R\to X}' } \ar@<.5ex>[r]^-{p_1} \ar@<-.5ex>[r]_-{p_2} &
{\Omega}_{U\to X}' }
\eqend
is an \'etale equivalence relation with quotient $X^n_S/\Sgot_n$.

Assume now that the base $S=\Spec(A)$ is affine, and that $X$ is a quasi-compact
algebraic space. Let $U=\Spec(F) \ra X$ be an \'etale affine cover. The map
$\Spec(\alpha_n) \colon \Spec(\TS^n_AF) \ra \Spec (\Gamma^n_AF)$ is a
(universal) homeomorphism (see e.g. \cite[Corollary 4.2.5]{rydh_fams}), and we
let
$$ \Omega_{U\to X}'' \subseteq \Spec(\Gamma^n_AF)$$
denote the open set given as the image of $\Omega_{U\to X}'$ by
$\Spec(\alpha_n)$. Our first aim is to prove that the homeomorphic image of
(\ref{eteqrelsym}) also forms an \'etale equivalence relation.

\begin{prop}\label{completion} 
Let $F\ra F'$ be an \'etale extension of $A$-algebras. Let $\xi \in \Omega_{F\to
F'}''$ be a point of $A$. Then the induced map of completions
$$
(\Gamma^n_AF)_{{f(\widehat{\xi})}} \ra (\Gamma^n_AF')_{\widehat{\xi}}
$$
is an isomorphism, where $f(\xi) $ is the image of $\xi$ by the induced map
$\Spec(\Gamma^n_AF') \ra \Spec (\Gamma^n_AF)$.
\end{prop}

\begin{proof} 
It suffices to show that there are ideals $I_1\subset \Gamma^n_AF$ and
$J_1\subset \Gamma^n_AF'$ contained in the  ideals corresponding to the
points $f(\xi)$ and $\xi$, respectively, such that the induced map of formal
neighborhoods \eqbeg \label{indmapcompl} \lim_{\longleftarrow}
(\Gamma^n_AF)/I_1^m \ra \lim_{\longleftarrow} (\Gamma^n_AF')/J_1^m
\eqend
is an isomorphism.  As the morphism $\Spec(\Gamma^n_AF') \ra \Spec(\TS^n_AF')$
is a homeomorphism the point $\xi$ lifts to a point of $\Spec (\T^n_AF')$. Let
$\varphi \colon \T^n_AF' \ra L'$ be a lifting of $\xi=\Spec(L)$, with $L'$ some
field extension of $L$. Write $\varphi =(\varphi_1, \ldots, \varphi_n)$, and
define the ideal $J=\cap \ker(\varphi_i)$ in $F$. We let $J_m =\ker
(\Gamma^n_AF' \ra \Gamma^n_A(F'/J^m)$. As the map $\Gamma^n_AF' \ra L$ factors
via $\Gamma^n_A(F/J)$ we have that $J_1$ is contained in the ideal
$\ker(\Gamma^n_AF' \ra L)$. We let $I_m =\ker(\Gamma^n_AF \ra
\Gamma^n_A(F/I^m))$ where $I=\cap \ker(\varphi_{i| F})$, and we consider the
induced map (\ref{indmapcompl}).

By Lemma (\ref{proeq}) we have the limit of the system $\{(\Gamma^n_AF)/I_1^m\}$
equals the limit of the system $\{ (\Gamma^n_AF)/I_m=\Gamma^n_A(F/I^m)\}$. By
Lemma (\ref{inftyngh}) we have that $F/I^m=F'/J^m$, and it follows that the map
(\ref{indmapcompl}) is an isomorphism.
\end{proof}

\begin{cor}\label{normext} 
Let $F \ra F'$ be an \'etale extension of $A$-algebras, and let $I_F\subseteq
\Gamma^n_AF$ and $I_{F'}\subseteq \Gamma^n_AF'$ be the ideal of norms associated
to $F$ and $F'$, respectively. These two ideals, $I_F\Gamma^n_AF'$ and $I_{F'}$,
are equal when restricted to the FPR-set $\Omega_{F\to F'}''\subseteq
\Spec(\Gamma^n_AF')$.
\end{cor}

\begin{proof} 
Assume first that the result is true when $F$ (and hence $F'$) is a
finitely presented $A$-algebra. We can write $f \colon F \ra F'$ as a
limit by a directed set of \'etale maps $f_{\alpha} \colon F_{\alpha}
\ra F'_{\alpha}$ of finitely presented $A$-algebras, such that
$F'_{\alpha}\otimes_{F_{\alpha}}F_{\beta}\backsimeq F'_{\beta}$ for
all $\alpha$ and all $\beta \geq \alpha$.  This means that
$\Spec(\T^n_AF') \to \Spec(\T^n_AF)$ can be thought of as
$\ili_\beta\Spec(\T^n_AF'_\beta) \to \ili_\beta\Spec(\T^n_AF_\beta)$
and similarly for $\T^n$ replaced by $\TS^n$ (as directed direct
limits commute with taking invariants) and $\Gamma^n$. The equality to
be proven is one of equality of stalks so we may focus on a particular
point $x'' \in \Omega ''_{F\to F'}$ which is the image of some $x \in
\Omega_{F\to F'}$. By Lemma (\ref{fpr-lemma}) we may assume that all
projection maps $\Spec(\T^n_AF) \to \Spec(\T^n_AF_\alpha)$ are \fpr{}
at the image of $x$ in $\Spec(\T^n_AF)$ and hence, again by Lemma
(\ref{fpr-lemma}), we get that $\Spec(\T^n_AF') \to \Spec(\T^n_AF)$ is
\fpr{} at $p_\alpha(x)$ for all $\alpha$ which means that
$p_\alpha(x'') \in \Omega ''_{F_\alpha \to F'_\alpha}$. Hence we get
that the equality $I_{F_\alpha}\Gamma^n_AF'=I_{F'_\alpha}\Gamma^n_AF'$
at $x''$ and taking the direct limit of sheaves in $\alpha$ gives the
Corollary at $x''$ and hence in $\Omega ''_{F\to F'}$.

We are therefore left with the case when $F$ is a finitely presented
$A$-algebra. By another (simpler) limit argument we reduce to the case when $A$
is Noetherian. To show equality of the two
ideals $I_F'$ and $I_F\Gamma^n_AF'$ it now suffices to show equality in the
completion $(\Gamma^n_AF')_{\widehat{\xi}}$, for each closed point $\xi \in
\Omega_{F\to F''}$ lying over the maximal ideal of $A$. The result now follows
from the proposition.
\end{proof}

\begin{cor}\label{etalness} 
Let $F\ra F'$ be an \'etale extension of $A$-algebras. The induced map
$\Omega_{F\to F'}'' \ra \Spec(\Gamma^n_AF)$ is \'etale.
\end{cor}

\begin{proof}  
By doing the same reductions as in the previous corollary, we may
assume that $A$ is Noetherian and $F$ and $F'$ are finite type
$A$-algebras. By localisation and Henselisation in $A$ we may then
also assume that $A$ is strictly Henselian. The result then follows
from the proposition.
\end{proof}

\begin{lem}\label{eteqrel} 
Let $X \ra S$ be a quasi-compact separated algebraic space over an
 affine base $S$. Write $X$ as a quotient $\xymatrix@M=1pt{R
 \ar@<.5ex>[r] \ar@<-.5ex>[r] & U}$, with affine schemes $U$ and
 $R$. Then we have that $\xymatrix@M=1pt{ \Omega_{R\to X}''
 \ar@<.5ex>[r] \ar@<-.5ex>[r] & \Omega_{U\to X}''}$ is an \'etale
 equivalence relation.
\end{lem}

\begin{proof} We have (\ref{eteqrelsym}) that $\xymatrix@M=1pt{ \Omega_{R\to X}' \ar@<.5ex>[r] \ar@<-.5ex>[r] & \Omega_{U\to X}'}$ is an \'etale equivalence relation. As the map $\Omega_{U\to X}' \ra \Omega_{U\to X}''$ is a homeomorphism, we have that the induced map $\Omega_{R\to X}'' \ra \Omega_{U\to X}''\times_S \Omega_{U\to X}''$ is injective over $\Spec(k)$-valued points, with $k$ a field. To prove the lemma it then suffices to show that the two maps 
$\xymatrix@M=1pt{ \Omega_{R\to X}'' \ar@<.5ex>[r] \ar@<-.5ex>[r] & \Omega_{U\to X}''}$
are \'etale. The projection maps $p_i \colon R \ra U$ are \'etale, and we have open immersions $\Omega_{R\to X} \subseteq \Omega_{p_i \colon R \to U}$. \'Etaleness of the two maps $\xymatrix@M=1pt{ \Omega_{R\to X}'' \ar@<.5ex>[r] \ar@<-.5ex>[r] & \Omega_{U\to X}''}$ then follows from Corollary (\ref{etalness}).
\end{proof}

\begin{prop}\label{cartesiandiagrams} Let $X\ra S$ be a separated quasi-compact algebraic space over an affine scheme $S=\Spec(A)$. Let $U=\Spec(F) \ra X$ be an \'etale affine cover, and let $R=U\times_XU$. Define $\Gamma^n_{X/S}$ as the quotient of the \'etale equivalence relation $ \xymatrix@M=1pt{ \Omega_{R\to X}'' \ar@<.5ex>[r] \ar@<-.5ex>[r] & \Omega_{U\to X}''}$. 
\begin{enumerate} 
\item We have a cartesian diagram
$$ \xymatrix{ 
\CalH^n_{U\to X} \ar[d]^{\n_U} \ar[r] & \Hilb^n_{U/S} \ar[d]^{\n_U} \\
\Omega_{U\to X}'' \ar[r] & \Gamma^n_{U/S}=\Spec(\Gamma^n_A(F)). }$$
\item In the diagram below we have $n_U\circ p_i =q_i\circ n_R \ ,i=1,2$, and consequently there is an induced map $\n_X \colon \Hilb^n_{X/S} \ra \Gamma^n_{X/S}$:
$$ 
\xymatrix{ {\CalH^n_{R\to X} } \ar[d]^{\n_R} \ar@<.5ex>[r]^-{p_1}
\ar@<-.5ex>[r]_-{p_2} & \CalH^n_{U\to X} \ar[d]^{\n_U} \ar[r]^-p &
\Hilb^n_{X/S}\ar[d]^{\n_X} \\ {\Omega_{R\to X}'' }\ar@<.5ex>[r]^-{q_1}
\ar@<-.5ex>[r]_-{q_2} & \Omega_{U\to X}'' \ar[r]^-{q} &
\Gamma^n_{X/S}}
$$
Moreover, the commutative diagrams above are cartesian.
\end{enumerate}

\end{prop}

\begin{proof} 
Let us first consider the special case with $S=\Spec(k)$, where $k$ is
an algebraically closed field. A $k$-valued point $Z\subseteq U$ of
the Hilbert functor $\Hilb^n_{U/S}$ has support at a finite number of
points $\xi_1, \ldots ,\xi_p$. By (\ref{sec1.10}) the associated cycle
$\n_U(Z)$ consist of the points $\xi_1, \ldots, \xi_p$ counted with
multiplicities $m_1, \ldots, m_p$. We have that the cycle $n_U(Z)$ is
in the \fpr-set $\Omega_{U\to X}''$ if and only if the closed
subscheme $Z\subseteq U$ also is a closed subscheme of $X$.

Now, let us prove the proposition. In the first diagram (1) the horizontal maps are open immersions. To see that it is commutative and cartesian it suffices to establish the equality of the two open sets $\CalH^n_{U\to X}$ and $\n_U^{-1}(\Omega_{U\to X}'')$ of $\Hilb^n_{U/S}$. This we can be checked by reducing to $S=\Spec(k)$, with $k$ algebraically closed. Then we are in the special case considered above from which Assertion (1) follows.

In particular we have proven that the restriction of the norm map $n_U$ to the  open subset $\CalH^n_{U\to X}$ has $\Omega_{U\to X}''$ as domain. We therefore obtain the two leftmost diagrams in (2).  Since the horizontal maps in these diagram are \'etale (Proposition (\ref{eqHilb}) and Lemma (\ref{eteqrel})) we can prove the diagrams are cartesian by evaluation over algebraically closed points. We are then again reduced to the special case considered above, which proves assertions in (2).
\end{proof}

\begin{prop}[Rydh]\label{rydh} Let $X \ra S$ be a separated map of algebraic spaces. Then there exists an algebraic space $\Gamma^n_{X/S} \ra S$ such that
\begin{enumerate}
\item When $X \ra S$ is quasi-compact with $S$ an affine scheme, the space
$\Gamma^n_{X/S}$ coincides with the one constructed above
(\ref{cartesiandiagrams}).  \item For any base change map $T\ra S$ we have a
natural identification $\Gamma^n_{X/S}\times_ST= \Gamma^n_{X\times_ST/T}$.
\item For any open immersion $X' \subseteq X$ we have an open immersion
$\Gamma^n_{X'/S}\subseteq \Gamma^n_{X/S}$, and moreover
$$ \Gamma^n_{X/S} = \lim_{
\underset{\text{open, q-compact}}{X'\subseteq X}}\Gamma^n_{X'/S}.$$ \item There
is a universal homeomorphism $X^n_S/\Sgot_n \ra \Gamma^n_{X/S}$, which is an
isomorphism when $X\ra S$ is flat, or when the characteristic is
zero.\end{enumerate}
\end{prop}

\begin{proof} 
All results can be found in (\cite{rydh_fams}): Existence of the space
$\Gamma^n_{X/S}$ is Theorem (3.4.1), whereas Assertion (4) is Corollary (4.2.5),
and the statement about open immersions in (3) is a special case of Proposition
(3.1.7). The functorial description of $\Gamma^n_{X/S}$ given by David Rydh
immediately gives assertion (2) and that $\Gamma^n_{X/S}$ is the union of
$\Gamma^n_{X'/S}$ with quasi-compact $X'\subseteq X$. Assertion (1) follows as
our $\Omega_{U\to X}''$ is what Rydh denotes with $\Gamma^n(U/S)_{|
\text{reg}/f}$ (see Proposition (4.2.4), and the proof of Theorem (3.4.1),
loc. cit.).
\end{proof}

\subsection{The ideal sheaf of norms} For $X\ra S$ quasi-compact and separated
over an affine base we have by Corollary (\ref{normext}) that the ideals of
norms patch together to form an ideal sheaf $\CalI_X$ on $\Gamma^n_{X/S}$. As
these ideals clearly commute with open immersions and base change we obtain by
(3) and (1) of Proposition (\ref{rydh}), an ideal sheaf of norms $\CalI_X$ on
$\Gamma^n_{X/S}$, for any separated algebraic space $X \ra S$. Let
$$ \Delta_X \subseteq \Gamma^n_{X/S}$$
denote the closed subspace defined by the ideal sheaf of norms.

\begin{thm} 
Let $X\ra S$ be a separated morphism of algebraic spaces. Then the good
component $\good$ of $\Hilb^n_{X/S}$ is isomorphic to the blow-up of
$\Gamma^n_{X/S}$ along the closed subspace $\Delta_X \subseteq \Gamma^n_{X/S}$,
defined by the ideal of norms associated to $X\ra S$. Moreover, if $X\ra S$ is
flat then $\good $ is obtained by blowing-up the geometric quotient
$X^n_S/\Sgot_n$.
\end{thm}

\begin{proof} 
The Hilbert scheme $\Hilb^n_{X/S}$ and $\Gamma^n_{X/S}$ commute with arbitrary
base change. The good component $\good$ as well as blow-ups, commute with flat,
and in particular \'etale base change. We may therefore assume that the base $S$
is an affine scheme.

For any open immersion $X'\subseteq X$, with $X'$ quasi-compact, we have a norm
map $\n_{X'} \colon \Hilb^n_{X'/S} \ra \Gamma^n_{X'/S}$ which, by varying $X'$,
form a norm map $\n_X \colon \Hilb^n_{X/S} \ra \Gamma^n_{X/S}$. We claim now
that the inverse image $\n^{-1}_X(\Delta_X)$ is locally principal, which we can
verify on an open cover. Moreover, given that we obtain an induced map from the
good component $\good$ to the blow-up of $\Gamma^n_{X/S}$ along $\Delta_X$. To
verify that the induced map is an isomorphism, we also reduce to an open
cover. Consequently we may assume that $X$ itself is quasi-compact.

When $X$ is quasi-compact we choose an \'etale affine cover $U\ra X$. Then by
using the cartesian diagrams (2) and (1) of Proposition
(\ref{cartesiandiagrams}) one establish using Theorem (\ref{thm3.6}) that
$\n_X^{-1}(\Delta_X)$ is locally principal. By Theorem (\ref{thm7.3}) we have
that the blow-up of $\Delta_U \subseteq \Gamma^n_{U/S}$ yields the good
component $\mathrm{G}^n_{U/S}$, and the isomorphism is given induced by the norm
map $\n_U$. It then follows by the two cartesian diagrams (2) and (1) of
Proposition (\ref{cartesiandiagrams}), that the map induced map from $\good$ to
the blow-up of $\Delta_X\subseteq \Gamma^n_{X/S}$ is an isomorphism.
\end{proof}

\subsection{The case of surfaces} 

Before we give a corollary to this result we need a generalisation of a result
of Fogarty on the smoothness of the Hilbert scheme (\cite[Theorem
2.9]{fogarty_algfam}). Fogarty proves that the Hilbert scheme of a smooth map $X
\ra S$ is smooth of relative dimension $2$ provided that $S$ is a Dedekind
scheme. As the Hilbert scheme commutes with base change and flatness can be
verified in the integral case by pulling back to Dedekind bases it follows that
the result of Fogarty is valid when the base $S$ is integral. However, as we
will see, no conditions on the base is needed for that statement. We shall give
a direct proof by proving formal smoothness using the infinitesimal lifting
criterion and the Hilbert-Burch theorem.

\begin{prop} Let $X \ra S$ be a smooth and separated morphism of
relative dimension $2$. Then ${\Hilb}^n_{X/S} \ra S$ is smooth for all $n$.
\end{prop}
\begin{proof} 
As $\Hilb^n_{X/S}$ commutes with base change we can assume that the base is
Noetherian. It is enough to show formal smoothness so the statement would follow
if we could show that for every small thickening $T \subset T'$ of local
Artinian $S$-schemes, any $T$-flat finite subscheme $Z \subseteq X\times_ST$ can
be extended to a $T'$-flat finite subscheme of $X\times_ST'$. Let $s$ be the
closed point in $S$. The obstruction for the existence of such a lifting is an
element $\alpha \in
\Ext^1_{\CalO_{X_s}}(\CalI_{Z_s},\CalO_{X_s}/\CalI_{Z_s})$. We have an exact
``local-to-global'' sequence
\begin{multline*}
H^1(X_s,\sHom_{\CalO_{X_s}}(\CalI_{Z_s},\CalO_{X_s}/\CalI_{Z_s})) \to
\Ext^1_{\CalO_{X_s}}(\CalI_{Z_s},\CalO_{X_s}/\CalI_{Z_s}) \to \\
H^0(X_s,\sExt^1_{\CalO_{X_s}}(\CalI_{Z_s},\CalO_{X_s}/\CalI_{Z_s})).
\end{multline*}

As $\sHom_{\CalO_{X_s}}(\CalI_{Z_s},\CalO_{X_s}/\CalI_{Z_s})$ has finite
support, the left term of the above sequence is $0$, and consequently it
suffices to show that the image of the obstruction element $\alpha$ in
$H^0(X_s,\sExt^1_{\CalO_{X_s}}(\CalI_{Z_s},\CalO_{X_s}/\CalI_{Z_s}))$ is
zero. As $Z$ is a disjoint union of points we have that $\alpha =\prod
\alpha_{z_i}$, where at a point $z \in Z$ the factor $\alpha_z$ is the
obstruction for lifting $\Spec \CalO_{Z,z}$, which is a closed flat subscheme of
$\Spec \CalO_{{X\times_ST},z}$, to a flat subscheme of $\Spec
\CalO_{{X\times_ST},z}$. It is thus enough to show that these local obstructions
vanish. Hence our situation is as follows: We have a surjection of local
Artinian rings $R' \ra R$ whose kernel is $1$-dimensional over the residue
field, an essentially smooth $2$-dimensional local $R'$-algebra $S'$, and a
quotient $S:=S'\bigotimes_{R'}R \ra T$ such that $T$ is a finite flat
$R$-module. We then want to lift $T$ to a quotient $S' \ra T'$ which is a flat
$R'$-module. We first claim that $T$ has projective dimension $2$ over $S$. As
$T$ is $R$-flat it is enough to check $\overline{T}$ has projective dimension
$2$ over $\overline{S}$, where $\overline{(-)}$ denotes reduction modulo the
maximal ideal of $R$. In that case we have that $\overline{T}$ is a
Cohen-Macaulay module over the regular local ring $\overline{S}$ with support of
codimension $2$ and the result follows.

By \cite[Thm.\ 7.15]{northcott} (cf. also the original proof in \cite{burch}) it then
follows that the ideal $I_T$ defining $T$ is the determinant ideal of $n\times
n$-minors of an $n+1\times n$-matrix $M$ and that the grade (the maximal length
of $S$-regular sequence contained in $I_T$) of $I_T$ is $2$. We then
(arbitrarily) lift $M$ to a matrix $M'$ over $S'$ and let $T'$ be defined by
$n\times n$-minors of $M'$. What remains to show is that $T'$ is $R'$-flat. The
grade of $I_{T'}$ is also $2$ as we may lift an $S$-regular sequence in $I_T$ to
elements of $I_{T'}$ which then given an $S'$-regular sequence and hence by
\cite[Thm.\ 7.16]{northcott}, the sequence
$$
0 \ra (S')^n \ra (S')^{n+1} \ra S' \ra T' \ra 0
$$
is exact, where $(S')^n \ra (S')^{n+1}$ is given by the lifted matrix and
$(S')^{n+1} \ra S'$ by its minors (with appropriate signs). For the same reason
this sequence tensored with the residue field of $R'$ remains exact which shows
that $T'$ is $R'$-flat. 
\end{proof}
\begin{cor} Let $X \ra S$ be a smooth, separated morphism of pure relative
dimension $2$. Then we have that the Hilbert scheme
${\Hilb}^n_{X/S}$ is the blow-up of $\Gamma^n_{X/S}$ along $\Delta_X$.
\end{cor}

\begin{proof} If we can prove that $U^{et}$ of ${\Hilb}^n_{X/S}$ is
schematically dense then we are finished by the Theorem. As the defining ideal
of the complement of $U^{et}$ is locally principal and as ${\Hilb}^n_{X/S} \ra
S$ is flat by the proposition this can be checked fibre by fibre so we may
assume that $S$ is the spectrum of a field $k$. Now, in that case
${\Hilb}^n_{X/S}$ is smooth by the proposition or by Fogarty's result. For the
density statement we may reduce to the base field $k$ being algebraically
closed. Write $X=\sqcup_{i=1,\ldots, p} X_i$ as a disjoint union. We have that
$\Hilb^n_{X/S}$ is the disjoint union $\sqcup_{n_1+\dots + n_p=n}
\Hilb^{n_i}(X_i)$. As $U^{et}$ is non-empty in each of the components
$\Hilb^{n_i}(X_i)$, this implies that it is schematically dense in
${\Hilb}^n_{X/S}$.
\end{proof}

\begin{rem} As pointed out by the referee, there is a small inaccuracy in
(\cite[Proposition 2.3]{fogarty_algfam}) concerning the connectedness of the
Hilbert scheme in that the Hilbert scheme of a connected scheme is
not necessarily connected. The proof had to take that into account.
\end{rem} 

\section{The good component for affine varieties}

We will in this last section generalize the approach Haiman gives in
\cite{haiman_catalan}, using the fact that the Hilbert scheme
$\operatorname{Hilb}^n_Y$, for a projective scheme $Y$, can be embedded as a
closed subscheme of the Grassmannian of rank $n$-quotients of
$H^0(Y,\CalO_Y(N))$, when $N$ is large enough. To simplify we assume that our
base $\Spec(A)$ is Noetherian.

\begin{prop}\label{finitephi} 
Let $X=\Spec(F)\ra S=\Spec(A)$ be a finite type morphism of affine schemes, and
let $V\subseteq F$ be an $n$-sufficiently big $A$-submodule. Let $I_V$ and $I_F$
be the ideals of norms associated to $V$ and $F$, respectively. The natural
morphism $\oplus_{m \geq 0}I_V^m \to \oplus_{m\geq 0}I_F^m$ induces a morphism
$$ 
\varphi \colon \good=\Proj(\oplus_{m\geq 0}I_F^m) \ra
\mathrm{Bl}_{I_V}(\Gamma^n_AF) =\Proj(\oplus_{m \geq 0}I_V^m)
$$
which is finite.
\end{prop}

\begin{proof}
Let $U$ resp.\ $U'$ be the complement of $\Spec(\Gamma^n_AF)$ in
$\Spec(\oplus_{m\geq 0}I_F^m)$ resp.\ $\Spec(\oplus_{m \geq
0}I_V^m)$. That the map on $\Proj$'s is well-defined means that the
map on spectra maps $U$ into $U'$. Assume therefore, by way of
contradiction, that we have a closed point $x$ of $U$ that does not
map into $U'$. This gives us a field valued point of
$\operatorname{Hilb}^n_{\Spec(F)/\Spec(A)}$, i.e., an $n$-dimensional
quotient $F\bigotimes_Ak \to R$. However, the assumption that the
image of $x$ does not lie in $U'$ means that the image of $V$ does not
span $R$. This however contradicts the assumption that $V$ is
$n$-sufficiently big.

For graded elements $f$ in a graded ring $R$ we let $D_+(f)$ denote
the basic open affine given as the spectrum of the degree zero part of
the localized ring $R_f$. We have, for any $f\in I_V$ that
$\varphi^{-1}(D_+(f))=D_+(f)$, hence the morphism $\varphi $ is an
affine morphism. Since $F$ is assumed of finite type it follows from
Lemma (\ref{modpolyring}) that $I_F$ is of finite type, and
consequently $\good$ is proper over $\Spec (\Gamma^n_AF)$. Since
$\mathrm{Bl}_{I_V}(\Gamma^n_AF)$ is separated it follows that $\varphi
$ is proper. Thus the morphism $\varphi$ is both proper and affine,
hence finite.
\end{proof}

When $V\subseteq F$ is $n$-sufficiently big we have an induced morphism 
$$ h \colon \Hilb^n_{X/S} \ra \mathrm{Grass}^n_V$$
from the Hilbert scheme to the Grassmannian.

\begin{lem}\label{philift} 
Let $X=\Spec(F) \ra S=\Spec(A)$ be of finite type, and let $V\subset F$ be
$n$-sufficiently big, finitely generated $A$-module. We have a commutative
diagram
$$ \xymatrix{
\good \ar[d]^{\varphi} \ar[r] & \Hilb^n_{X/S} \ar[d]^{h}\\ \mathrm
{Bl}_{I_V}(\Gamma^n_AF) \ar[r] & \mathrm{Grass}^n_V}.
$$
\end{lem}

\begin{proof}  
Since $V$ is finitely generated we can use the Pl\"ucker coordinates to embed
$\mathrm{Grass}^n_V$ as a closed subscheme of ${\mathbf P}(\wedge^n
V)$. Composed with the diagonal embedding and the Segre embedding yields the
closed immersion $\iota_1$ given as the composite
$$ 
\mathrm{Grass}^n_V \subset {\mathbf P}(\wedge^nV) \subset {\mathbf P}(\wedge
^nV)\times {\mathbf P}(\wedge^nV) \subset {\mathbf P}(\wedge^nV \otimes
\wedge^nV).
$$
The natural map of $A$-modules $\wedge^nV\otimes_A\wedge^nV\ra I_V$ will by
definition hit all the generators for the ideal $I_V$, and consequently
determine a closed immersion $\iota_2 \colon \mathrm{Bl}_{I_V}(\Gamma^n_AF) \ra
\mathbf{P}(\wedge^n_V \otimes \wedge ^nV)\times\Spec(\Gamma^n_A(F))$. We now
have the commutative diagram
$$\xymatrix{
\good \ar[d]^{\varphi} \ar[r] & \Hilb^n_{X/S} \ar[r]^{h} & \mathrm{Grass}^n_V \ar[d]^{\iota_1}\\
\mathrm {Bl}_{I_V}(\Gamma^n_AF) \ar[rr]^{p_1\circ\iota_2} & & {\mathbf P}(\wedge^nV\otimes_A \wedge^nV)},
$$
where $p_1$ is the projection on the first factor. The inverse image
$\varphi^{-1}(E)$ of the exceptional divisor $E\subseteq
\mathrm{Bl}_{I_V}(\Gamma^n_AF)$ is the exceptional divisor of $\good$, and on
the open complement we have that $\varphi $ is an isomorphism. Consequently
$p_1\circ\iota_2 \colon \mathrm{Bl}_{I_V}(\Gamma^n_AF) \ra {\mathbf P}(\wedge^nV
\otimes_A \wedge^nV)$ factors through $\mathrm{Grass}^n_V$ since it does so on
the complement of a Cartier divisor.
\end{proof}

\subsection{} 
Consider now $Y={\bf P}^r_S$, and let $g \colon Y \ra S$ denote the structure
map. For any closed subscheme $Z\subseteq Y$ that is flat, locally free of rank
$n$ over $S$, the induced map \eqbeg \label{suffbig} g_*\CalO_Y(N) \ra
g_*\CalO_Z(N)
\eqend
is easily seen to be surjective for $N\geq n-1$. Furthermore, the ideal sheaf
$\CalI_Z$ twisted with $N\geq n$ is regular, that is $R^pg_*\CalI_Z(N-p)=0$
for $p>0$ when $N\geq n$. It follows (\cite{gr_existence}) that the induced
morphism
\begin{equation}\label{Hilbimmersion}
\Hilb_{Y/S}^n \ra \mathrm{Grass}^n_{g_*\CalO_Y(N)}
\end{equation}
is a closed immersion for $N\geq n$.

\begin{prop} 
Let $F$ be an $A$-algebra generated by $t_1, \ldots, t_r$, let $V\subseteq F$ be
spanned by the monomials of degree $\leq n$ in the $t_1, \ldots, t_r$. Then the
morphism
$$ \varphi \colon \good \ra \mathrm{Bl}_{I_V}(\Gamma^n_AF)$$
is an isomorphism.
\end{prop}

\begin{proof} 
We embed $X=\Spec (F)$ in $Y={\bf P}^r_S$ using $(1\colon t_1\colon \cdots
\colon t_r)$. The natural map $h \colon \Hilb^n_{X/S} \ra \mathrm{Grass}^n_V$ is
an immersion and the natural map
$\mathrm{Grass}^n_V\to\mathrm{Grass}^n_{g_{*}(\CalO_Y(N))}$ is a closed
immersion. As $\Hilb^n_{X/S}$ immerses into $\Hilb^n_{Y/S}$, and the map
(\ref{Hilbimmersion}) is an immersion, it follows that the map $h \colon
\Hilb^n_{X/S} \ra \mathrm{Grass}^n_V$ is an immersion.

By Lemma (\ref{philift}) we  have that the restriction of $h$ to $\good $ factors through 
$$ \varphi \colon \good \ra \mathrm{Bl}_{I_V}(\Gamma^n_AF),$$
hence $\varphi $ must be an immersion as well. However, by Proposition
(\ref{finitephi}) the map $\varphi$ is proper, and consequently we have that the
map $\varphi$ must be a closed immersion. Furthermore, since $\varphi $ is an
isomorphism over the complement of a Cartier divisor, it is an isomorphism.
\end{proof}


\bibliographystyle{plain}
\bibliography{space_gen_et_fam}

\end{document}